\documentclass[letterpaper]{article}

\usepackage{arxiv}
\usepackage[    
    backend=biber,
    style=numeric-comp,
    sorting=none,
    giveninits=true,
    maxbibnames=99,
    maxcitenames=99
    ]{biblatex}
\bibliography{bibliography}
\usepackage{amsmath, amssymb, amsthm, amsfonts, graphicx, lineno, multicol, hyperref, xcolor, mathtools}
\usepackage[normalem]{ulem}
\usepackage{caption} 
\usepackage{subcaption}

\usepackage{macros}

\hypersetup{
pdftitle={PolyChopper: a Polyhedron Splitting Scheme},
pdfsubject={},
pdfauthor={Tommaso~Sorgente, Fabio~Vicini},
pdfkeywords={},
}

\begin{document}

\title{PolyChopper: a Polyhedron Splitting Scheme}

\author{
  \href{https://orcid.org/0000-0001-8699-6379}{\includegraphics[scale=0.06]{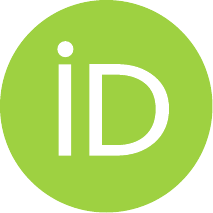}\hspace{1mm}Tommaso~Sorgente} \\
	Istituto di Matematica Applicata e Tecnologie Informatiche\\
	``E. Magenes''\\
	CNR, Genova, 16149, Italy \\
	\texttt{tommaso.sorgente@cnr.it} \\
  \And
	\href{https://orcid.org/0000-0001-7123-9199}{\includegraphics[scale=0.06]{orcid.pdf}\hspace{1mm}Fabio~Vicini} \\
	Dipartimento di Scienze Matematiche\\
	``G. L. Lagrange''\\
	Politecnico di Torino, Torino, 10129, Italy \\
	\texttt{fabio.vicini@polito.it} \\
}
\maketitle

\begin{abstract}
We introduce and implement a novel scheme, called \textit{PolyChopper}, for splitting a convex polyhedron into a finite set of disjoint convex sub-polyhedra.
The algorithm cuts the polyhedron with a plane, which can either be prescribed as an external constraint or automatically determined from the inertia tensor.
The resulting intersection vertices are then adjusted to control the complexity and quality of the generated elements, introducing additional ``notch wedges'', tetrahedra and pyramids, according to a user-defined quality parameter.
Experimental results demonstrate that the proposed approach is robust and consistently produces sub-polyhedra that are both simpler and of higher quality than the original polyhedron.
Furthermore, recursively applying the algorithm generates a hierarchy of polyhedral subdivisions with progressively smaller elements, all guaranteed to remain convex and predominantly tetrahedral.
\end{abstract}

\keywords{polyhedron splitting; polygonal mesh; polyhedral mesh; mesh generation; mesh refinement.}

\section{Introduction}
\label{sec:intro}
We consider the problem of splitting a convex polyhedron into a finite number of disjoint convex parts.
This problem, commonly referred to as a ``subdivision scheme,'' can be applied recursively in both local and global mesh generation and refinement processes, where the goal is to progressively reduce the size of the elements in a discrete tessellation starting from an initial polyhedron.

\begin{figure}[!ht]
    \centering
    \begin{subfigure}{0.32\linewidth}
        \centering
        \includegraphics[width=0.6\linewidth]{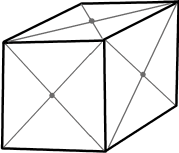}
        \caption{Tetrahedrization with inner points.}
        \label{fig:intro:existing:a}
    \end{subfigure}
    \begin{subfigure}{0.32\linewidth}
        \centering
        \includegraphics[width=0.65\linewidth]{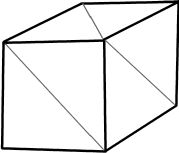}
        \caption{Tetrahedrization with ear-clipping.}
        \label{fig:intro:existing:b}
    \end{subfigure}
    \begin{subfigure}{0.32\linewidth}
        \centering
        \includegraphics[width=0.9\linewidth]{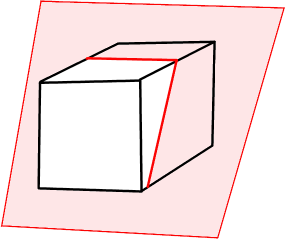}
        \caption{Brutal cut constrained to a plane.}
        \label{fig:intro:existing:b}
    \end{subfigure}
    \caption{Existing approaches for splitting a convex polyhedron.}
    \label{fig:intro:existing}
\end{figure}

Existing approaches to this problem include the simple free-splitting scheme known as \emph{barycentric subdivision}~\cite{comp_top}, in which the barycenter of every geometric object is computed and connected according to the face lattice to produce a subdivision into sub-polyhedra.
A commonly used extension first triangulates each face and then decomposes the polyhedron into tetrahedra by connecting the polyhedron barycenter to the resulting face triangles; see Figures~\ref{fig:intro:existing:a} and~\ref{fig:intro:existing:b}.
This operation can be performed either using an ear-clipping scheme, which is more robust, or by selecting an interior point, e.g., the centroid of each face, and connecting it to the vertices, which generally produces higher-quality elements.

Other classical approaches widely used in applications are \emph{dual-based} subdivision methods, such as Voronoi diagrams~\cite{voro} and Delaunay tessellations~\cite{tetgen}.
In recent years, an other algorithm belonging to this class have been proposed based on novel neural network techniques; see, for example,~\cite{ANTONIETTI2022111531}.

Finally, we mention the class of constrained subdivision schemes known as \emph{Binary Space Partitions} (BSP)~\cite{deBerg2008}, in which the subdivision of the polyhedron is guided by prescribed directions, e.g., a set of planes that split the original volume into subdomains.
This class includes what we refer to throughout the paper as a \textit{brutal cut}, namely, cutting the polyhedron with a pre-selected plane without considering the quality of the resulting objects; see Figure~\ref{fig:intro:existing:b}.
The cutting plane may be prescribed as a constraint or selected to optimize a quality criterion, for instance by exploiting the inertia tensor of polygonal mesh elements~\cite{berrone2025effective}.
Although this approach is fast and simple, it often generates tiny or low-quality polyhedra, faces, or edges, making the resulting subdivision unsuitable for practical applications.
Moreover, particularly complex intersections, e.g., when the cutting plane passes very close to a vertex, may give rise to numerical issues that compromise the robustness of the algorithm.

In this work, we introduce a novel subdivision scheme, named \emph{PolyChopper}, based on the BSP class, that extends the brutal cut approach with the following objectives:
\begin{itemize}
    \item preserving the quality of the resulting sub-polyhedra, as measured by a selected quality indicator, see \cite{SBMS23};
    \item minimizing the \emph{complexity} of the generated sub-polyhedra, according to a suitable definition of complexity.
\end{itemize}

The name \emph{PolyChopper} is inspired by the technique used by woodcutters to cut (chop) large trees into manageable pieces. Similarly, the proposed strategy subdivides a polyhedron into two ``main'' high-quality convex polyhedra lying on opposite sides of the cutting plane by introducing a sequence of small ``notches'' along the plane. These notches remove a collection of small ``notch wedges'' (tetrahedra and pyramids) located between the two main polyhedra.

The wedges are generated according to a quality indicator for polyhedral elements and are designed to preserve the convexity of the two main polyhedra. Moreover, the purpose of the notches is to eliminate or relocate the intersection points away from the vertices of the original polyhedron, thereby preventing poorly shaped elements that may lead to numerical instability, reducing the number of newly introduced vertices, and improving the quality of the two main polyhedra.
Indeed, an important aspect of any subdivision scheme is the number of new objects (polyhedra, faces, edges, and vertices) generated during the subdivision process. Ideally, a polyhedron should be partitioned into the smallest possible number of convex polyhedra while maintaining a high level of quality for all the newly generated objects.

In summary, the main features of the proposed algorithm are as follows:
\begin{itemize}
    \item \textit{robustness:} the introduction of the wedges avoids complex intersections that may lead to numerical instability;
    \item \textit{efficiency:} the algorithm achieves a good balance between the number of generated objects and their quality, while avoiding pathological configurations;
    \item \textit{flexibility:} the algorithm can be driven by any polyhedral quality indicator and naturally reduces to the brutal cut algorithm when no quality indicator is specified.
\end{itemize}

The proposed algorithm extends the two-dimensional approach introduced in \cite{BERRONE2021103502} and further developed in \cite{berrone2025effective}. In those works, the splitting line is adjusted within the polygon to improve element quality and reduce the number of vertices. However, the three-dimensional problem, here tackled, is significantly more challenging. While a line intersects a convex polygon at most at two points, a plane may intersect a convex polyhedron along a polygon with more than three vertices. Consequently, modifying the intersection points to improve the quality of the resulting elements preserves a valid splitting line in 2D, but generally does not define a valid splitting plane in 3D, as we will shown in the following sections.

We are confident that the proposed constrained splitting scheme will prove valuable for coupled problems involving highly complex geometries. In particular, in our future research it will improve the three-dimensional adaptation of the tessellation to low-dimensional features of the domain, such as fractures in porous media \cite{10.1007/s10092-023-00517-5}, blood vessels in the human body \cite{berrone2023optimization}, and tree roots in the soil \cite{berrone20263d}.

Finally, we would like to remark that we are aware that requiring the polyhedra to remain convex may be considered a strong constraint, as preserving convexity throughout the splitting process limits the number of admissible operations.
However, we believe that maintaining convexity during the splitting process is crucial for numerical simulations, as it makes both the geometric algorithms and the discretization methods significantly more robust. Indeed, for convex polyhedral elements, numerical integration is well defined, and classical numerical methods, such as the Finite Element Method (FEM) and the Finite Volume Method (FVM), generally lead to better-conditioned linear systems.

\paragraph{Notation}
Given a $\dimension$-dimensional polytope $E$, namely a linear segment if $\dimension=1$, a planar polygon if $\dimension=2$, a polyhedron if $\dimension=3$, we denote $\#E$ the number of its $(\dimension-1)$-dimensional polytopes, i.e. vertices if $\dimension=1$, edges if $\dimension=2$, faces if $\dimension=3$.
We then indicate with $\|E\|$ its $\dimension-$dimensional size (length, area, or volume), with $|E|$ its $(\dimension-1)$-dimensional size (perimeter or surface area), and with $h_E$ its diameter, defined as the maximum distance between two vertices of $E$.
Throughout the paper, we will use the symbols $\polyhedron$ and $\polygon$ for polyhedra and polygons, and $\face$, $\edge$, and $\vert$ for faces, edges, and vertices, respectively.

\paragraph{Organization}
The remainder of the paper is organized as follows.
In Section~\ref{sec:algorithm} we describe the algorithm in detail, through visual interpretations.
In Section~\ref{sec:results} we test the algorithm in different scenarios, proving its robustness and efficiency.
Section~\ref{sec:conclusions} summarizes our findings, discusses their implications, and outlines potential directions for future research.

\section{The PolyChopper Algorithm}
\label{sec:algorithm}
Given a convex polyhedron $\polyhedron\subset\mathbb{R}^3$ and a plane $\plane$ intersecting $\polyhedron$, the PolyChopper algorithm splits $\polyhedron$ along the cut defined by $\plane$ into a set of convex sub-polyhedra $\{\polyhedron_1, \ldots, \polyhedron_s\}$, where $s \geq 2$ is finite. The algorithm is designed to control the splitting process so as to preserve the quality of the resulting sub-polyhedra, as measured by a selected quality indicator from the ones proposed in \cite{SBMS23}, ensuring, if possible, that it remains as high as, or comparable to, that of the original polyhedron. Furthermore, the algorithm seeks to minimize the complexity of the generated sub-polyhedra. We define the \emph{complexity} of a polytope $E$ according to its geometric dimension $\dimension$. For a polygon $(\dimension = 2)$, complexity is uniquely determined by the number of vertices; thus, polygons with fewer vertices are considered less complex.
For a polyhedron $(\dimension = 3)$, we consider both the number of vertices and the number of faces. Since the relative importance of these two quantities depends on the context, we regard one polyhedron as simpler than another only if it has both fewer vertices and fewer faces.

The intersection between $\polyhedron$ and $\plane$ defines a convex planar polygon $\polygon_0 \coloneqq \polyhedron\cap\plane=\left\{\vert_1,\ldots,\vert_n\right\}$, where each vertex $\vert_i\in\polygon_0$ either lies on an edge $\edge_j$ of $\polyhedron$ or coincides with a vertex $\vert_k$ of $\polyhedron$.
At each vertex $\vert_i \in \polygon_0$, we assign a distance $\distance_i \in \mathbb{R}$, denoting $\distanceset=\left\{\distance_1,\ldots,\distance_n\right\}$ the collection of these distances.
If $\vert_i$ coincides with a vertex of $\polyhedron$, we set $\distance_i = 0$; otherwise, we define a transformation $\transformation_i:\mathbb{R}^{3}\to\mathbb{R}^{3}$ that shifts $\vert_i$ along the edge $\edge_j$ on which $\vert_i$ lies, by a distance $\distance_i>0$.
Each $\distance_i$ is chosen to be sufficiently small to ensure that the transformed point remains on the closure of the corresponding edge $\bar{\edge}_j$.
With a slight abuse of notation, we denote the transformed point $\vert_i+\distance_i\edge_j$ simply by $\vert_i+\distance_i$.
Since the distances $\distance_i$ may vary independently for $i\in\{1,\ldots,n\}$, the resulting polygon $\polygon_{\transformation}\coloneqq\transformation(\polygon_0)$ is, in general, no longer planar.
Nevertheless, because each vertex is constrained to move only along its corresponding edge, the transformation cannot introduce self-intersections in $\polygon_{\transformation}$.

As stated at the beginning of this section, the distances $\distanceset$ are computed to achieve two objectives. The primary objective is to ensure that the newly generated sub-polyhedra have a quality, as measured by the prescribed quality indicator, that is comparable to or better than that of the original polyhedron $\polyhedron$. The second objective is to reduce the complexity of the sub-polyhedra introduced by the splitting process. To this end, each $\distance_i$ is chosen so that the shifted vertex $\vert_i+\distance_i$ coincides, whenever possible, with one of the endpoints of $\edge_j$ or with its barycenter.

In the following, for ease of presentation of the algorithm, we focus on the specific case illustrated in Figure~\ref{fig:algorithm:single_chop:a}, where the plane $\plane$ intersects $\polyhedron$ at $n = 5$ points along its edges. We denote these intersection points set by $\{\vert_{i-2},\vert_{i-1},\vert_i,\vert_{i+1},\vert_{i+2}\}$, the corresponding edges set by $\{\edge_{i-2},\edge_{i-1},\edge_i,\edge_{i+1},\edge_{i+2}\}$, and we analyze all possible configurations of the distance vector $\distanceset$.


\subsection{No Chop}
\label{sec:algorithm:none}
In the trivial case where all distances are zero, i.e., $\distance_i=0, \forall i \in \{1, \dots, n\}$, the polygon $\polygon_{\transformation}$ coincides with the original $\polygon_0$, and the plane $\plane$ splits $\polyhedron$ into the set $\{\polyhedron_1, \polyhedron_2\}$, that contains two convex polyhedra, see Figure~\ref{fig:algorithm:single_chop:a}.
In particular, we denote by $\polyhedron_1$ the polyhedron lying on the side of $\plane$ pointed to by the plane normal, shown on the right in the figure, and by $\polyhedron_2$ the polyhedron on the opposite side.

\subsection{Single Chop}
\label{sec:algorithm:single}
\begin{figure}[!ht]
    \centering
    \begin{subfigure}{0.45\linewidth}
        \centering
        \includegraphics[width=0.9\linewidth]{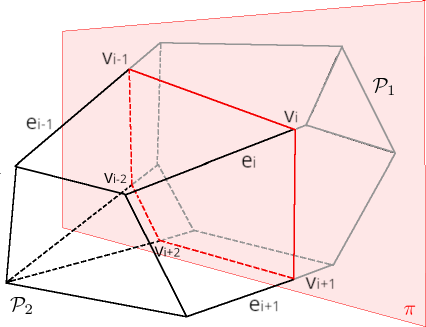}
        \caption{}
        \label{fig:algorithm:single_chop:a}
    \end{subfigure}
    \begin{subfigure}{0.45\linewidth}
        \centering
        \includegraphics[width=0.9\linewidth]{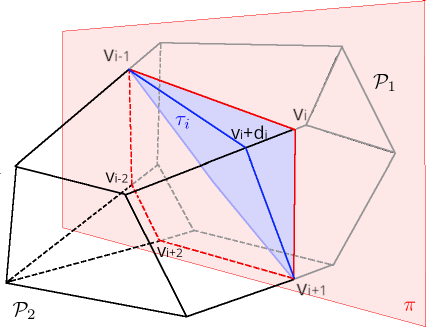}
        \caption{}
        \label{fig:algorithm:single_chop:b}
    \end{subfigure}
    \caption{(a) Simple polyhedron $\polyhedron$ split by the plane $\plane$. (b) Single chopping of $\polyhedron$ by edge $\edge_i$.}
    \label{fig:algorithm:single_chop}
\end{figure}

If only one distance is non-zero, e.g., $\distance_i\ne0$, the point $\vert_i$ is shifted along the edge $\edge_j$ to the new position $\vert_i+\distance_i$, as illustrated in Figure~\ref{fig:algorithm:single_chop:b}, which shows the case $\distance_i < 0$. The polygon $\polygon_{\transformation}\coloneqq\{\vert_{i-2},\vert_{i-1},\vert_i+\distance_i,\vert_{i+1},\vert_{i+2}\}$ is no longer planar. 
To restore convexity, we split $\polyhedron$ into three convex polyhedra: $\polyhedron_1$, containing $\vert_i$; $\polyhedron_2$, containing $\vert_i+\distance_i$; and $\polyhedron_3$, a ``chopped'' tetrahedron with vertices $\tet_i \coloneqq \{\vert_{i-1}, \vert_i,\vert_i+\distance_i,\vert_{i+1}\}$.

We emphasize that the common face between $\polyhedron_1$ and $\polyhedron_2$, which lies on $\plane$, has its complexity reduced by one vertex. Specifically, the original polygon $\polygon_0$ is replaced by the new polygon $\face \coloneqq \{\vert_{i-2},\vert_{i-1},\vert_{i+1},\vert_{i+2}\}$. This simplification is directly reflected in the complexity of $\polyhedron_1$ and $\polyhedron_2$. In particular, the new faces of $\polyhedron_1$ are the polygon $\face$ and the triangle $\{\vert_{i-1},\vert_{i},\vert_{i+1}\}$, while the new faces of $\polyhedron_2$ are the polygon $\face$ and the triangle $\{\vert_{i-1},\vert_{i} + \distance_i,\vert_{i+1}\}$.

\subsection{Double Chop}
\label{sec:algorithm:double}
\begin{figure}[!ht]
    \centering
    \begin{subfigure}{0.45\linewidth}
        \centering
        \includegraphics[width=0.9\linewidth]{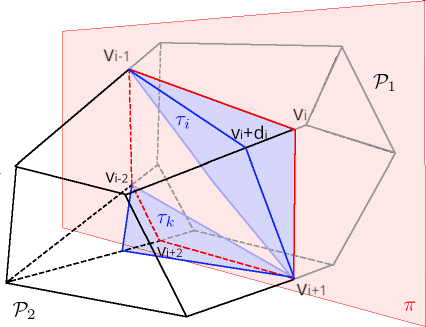}
        \caption{Non-consecutive $i,k =i+2$.}
        \label{fig:algorithm:double_chop:a}
    \end{subfigure}
    \begin{subfigure}{0.45\linewidth}
        \centering
        \includegraphics[width=0.9\linewidth]{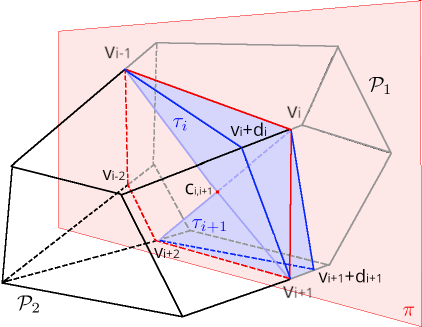}
        \caption{Consecutive $i,k =i+1$ with $\distance_i \distance_{i+1}<0$.}
        \label{fig:algorithm:double_chop:b}
    \end{subfigure}
    \begin{subfigure}{0.48\linewidth}
        \centering
        \includegraphics[width=0.9\linewidth]{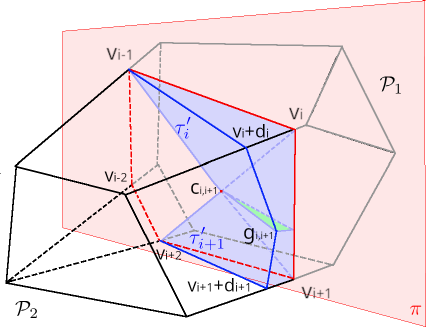}
        \caption{Consecutive $i,k = i+1$ with $\distance_i \distance_{i+1}>0$: straightforward application of the double chop.}
        \label{fig:algorithm:double_chop:c}
    \end{subfigure}
    \begin{subfigure}{0.48\linewidth}
        \centering
        \includegraphics[width=0.9\linewidth]{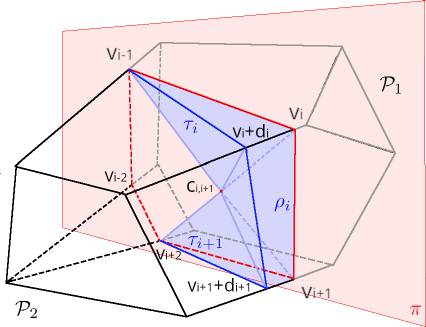}
        \caption{Consecutive $i,k =i+i$ with $\distance_i \distance_{i+1}>0$.}
        \label{fig:algorithm:double_chop:d}
    \end{subfigure}
    \caption{Double chopping of $\polyhedron$ by edges $\edge_{i}$ and $\edge_{k}$.}
    \label{fig:algorithm:double_chop}
\end{figure}

We now consider the case in which two distances are non-zero, namely $\distance_i\ne0$ and $\distance_k\ne0$, with $i\ne k$.
Figures~\ref{fig:algorithm:double_chop:a}--\ref{fig:algorithm:double_chop:c} illustrate all possible configurations.
Applying the single-chop construction of Section~\ref{sec:algorithm:single} to $\polyhedron$ along edge $\edge_i$ produces two convex polyhedra, $\polyhedron_1'$ and $\polyhedron_2'$, together with the tetrahedron $\tet_i$. The single-chop construction can then be applied along edge $\edge_k$ to either $\polyhedron_1'$ or $\polyhedron_2'$, depending on the sign of $\distance_k$. This generates the tetrahedron $\tet_k$ and the final convex polyhedra: if $\distance_k > 0$, then $\polyhedron_1=\polyhedron_1'-\tet_k$ and $\polyhedron_2=\polyhedron_2'$; otherwise, $\polyhedron_1=\polyhedron_1'$ and $\polyhedron_2=\polyhedron_2'-\tet_k$.
By construction, the two polyhedra $\polyhedron_1$ and $\polyhedron_2$ are convex and disjoint, since they are obtained through successive applications of the single-chop construction. In contrast, the two tetrahedra $\tet_i$ and $\tet_k$ may intersect, depending on the relative positions of the shifted vertices. Furthermore, the face lying on $\plane$ and shared by $\polyhedron_1$ and $\polyhedron_2$ may vary according to the values of $\distance_i$ and $\distance_k$.

If $k\neq \{i-1,i+1\}$, i.e., the two diagonals are not consecutive, the tetrahedra $\tet_i$ and $\tet_k$ cannot intersect because $\polygon_0$ is convex.
The result of the double chop is the set of four convex polyhedra $\{\tet_i,\tet_k,\polyhedron_1,\polyhedron_2\}$. Figure~\ref{fig:algorithm:double_chop:a} illustrates the case $k=i+2$.

When the two significant distances are consecutive and have opposite signs, e.g., $\distance_i \distance_{i+1}<0$, the transformations $\transformation_i$ and $\transformation_{i+1}$ shift the vertices $\vert_i$ and $\vert_{i+1}$, respectively, to opposite sides of the plane $\plane$, see Figure~\ref{fig:algorithm:double_chop:b}.
In this case, the two tetrahedra do not intersect because they lie on opposite sides of $\plane$. However, the edges $(\vert_{i+1},\vert_{i-1})$ and $(\vert_{i+2},\vert_{i})$ intersect at a point $\pcenter^i_{i+1}\in\plane$.
On the other hand, if for example $\distance_i \distance_{i+1}>0$, the transformations $\transformation_i$ and $\transformation_{i+1}$ shift $\vert_i$ and $\vert_{i+1}$, respectively, to the same side of $\plane$.
In this case, the intersecting tetrahedra $\tet_i$ and $\tet_{i+1}$ form a concave polyhedron, as shown in Figure~\ref{fig:algorithm:double_chop:c}.
In the figure, the vertex at which the concavity occurs is denoted by $\mathbf{g}^{i+1}_i$.
Note that one of the vertices of this polyhedron is the same point $\pcenter^{i+1}_i$ obtained in the case $\distance_i \distance_{i+1}<0$.

We can resolve both of these complex configurations by introducing two new edges connecting $\pcenter^{i+1}_i$ to the vertices $\vert_i+\distance_i$ and $\vert_{i+1}+\distance_{i+1}$, respectively; see Figure~\ref{fig:algorithm:double_chop:d} for the case $\distance_i\distance_{i+1} > 0$.
The case $\distance_i\distance_{i+1} < 0$ can be handled similarly.With these new edges, the resulting decomposition consists of two tetrahedra $\tet^{i}_{i-1}$ and $\tet^{i+2}_{i+1}$, and one pyramid $\pir^{i+1}_i$ with a quadrilateral base. More precisely:

\begin{align*}
    \tet^{i}_{i-1}   & \coloneqq \{\vert_{i-1},\vert_i,\vert_i+\distance_i,\pcenter^{i+1}_i\},                             \\
    \pir^{i+1}_i     & \coloneqq \{\vert_i,\vert_i+\distance_i,\vert_{i+1},\vert_{i+1}+\distance_{i+1},\pcenter^{i+1}_i\}, \\
    \tet^{i+2}_{i+1} & \coloneqq \{\vert_{i+1},\vert_{i+1}+\distance_{i+1},\vert_{i+2},\pcenter^{i+1}_i\}.
\end{align*}

Note that, the difference between this new proposal and the \emph{straightforward} double chop is the tetrahedron:

\begin{equation*}
    (\tet^i_{i-1}+\pir^{i+1}_i+\tet^{i+2}_{i+1}) - (\tet_i+\tet_{i+1}) = (\pcenter^{i+1}_i, \vert_i+\distance_i, \mathbf{g}^{i+1}_i, \vert_{i+1}+\distance_{i+1}),
\end{equation*}
which can be subtracted from $\polyhedron_2$ without compromising its convexity.
To summarize, when the two shifted vertices are consecutive, the resulting decomposition consists of five polyhedra, namely
$\{\tet^i_{i-1}, \pir^{i+1}_i, \tet^{i+2}_{i+1},\polyhedron_1, \polyhedron_2\}$.
We remark that all the new faces introduced outside the plane by the insertion of the new polyhedra are either triangles or quadrilaterals.
Moreover, the common face $\face$ shared by the two polyhedra $\polyhedron_1$ and $\polyhedron_2$ is simplified by the removal of one vertex, since two original vertices are removed and replaced by the new central vertex $\pcenter^{i+1}_i$.

\subsection{Triple Chop}
\label{sec:algorithm:triple}
\begin{figure}[!ht]
    \centering
    \begin{subfigure}{0.45\linewidth}
        \centering
        \includegraphics[width=0.9\linewidth]{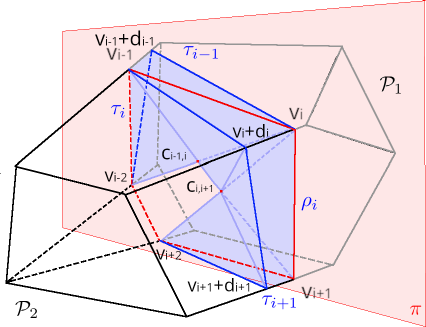}
        \caption{$\distance_{i-1}\distance_i<0$ and $\distance_{i}\distance_{i+1}>0$.}
        \label{fig:algorithm:triple_chop:a}
    \end{subfigure}
    \begin{subfigure}{0.45\linewidth}
        \centering
        \includegraphics[width=0.9\linewidth]{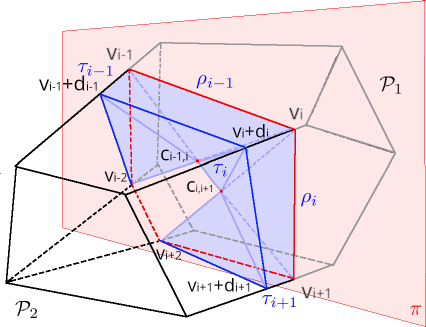}
        \caption{$\distance_{i-1}\distance_i>0$ and $\distance_{i}\distance_{i+1}>0$.}
        \label{fig:algorithm:triple_chop:b}
    \end{subfigure}
    \caption{Chopping of $\polyhedron$ with three consecutive edges $\edge_{i-1},\edge_i,\edge_{i+1}$.}
    \label{fig:algorithm:triple_chop}
\end{figure}

We discuss the case in which three distances are different from zero.
We omit the description of shifted vertices that are not consecutive, as their treatment is straightforward.
Instead, Figure~\ref{fig:algorithm:triple_chop} illustrates the case of three consecutive nonzero values $\distance_{i-1}, \distance_i, \distance_{i+1}$.
In both cases shown in the figure, namely when $\distance_{i-1}\distance_i<0$ or $\distance_{i-1}\distance_i>0$, we propose a solution analogous to that presented in the double-chop Section~\ref{sec:algorithm:double}, generalized to handle three intersections.
Accordingly, as depicted in Figure~\ref{fig:algorithm:triple_chop:b}, the algorithm produces two convex polyhedra $\polyhedron_1,\polyhedron_2$, two tetrahedra $\tet^{i-1}_{i-2}$ and $\tet^{i+2}_{i+1}$, two pyramids $\pir^i_{i-1}$ and $\pir^{i+1}_i$, and a new tetrahedron between them:

\begin{equation*}
    \tet^c_i=\left\{\pcenter^{i+1}_i,\pcenter^i_{i-1},\vert_i,\vert_i+\distance_i\right\}.
\end{equation*}

Overall, the resulting decomposition consists of seven sub-polyhedra, namely
$\{\tet^{i-1}_{i-2}, \pir^i_{i-1}, \tet^c_i, \pir^{i+1}_i, \tet^{i+2}_{i+1},\polyhedron_1, \polyhedron_2\}$.
Moreover, the common face $\face$ shared by the two polyhedra $\polyhedron_1$ and $\polyhedron_2$ is simplified by the removal of a single vertex, as three original vertices are removed and replaced by two new central vertices.

\subsection{General Multiple Chop}
\label{sec:algorithm:multiple}
The last case we need to discuss is when we have a sequence of $n$ nonzero entries $\distance_i,\ldots,\distance_{i+n}$.
In this case, generalizing the considerations reported till now,  the PolyChopper algorithm splits the input polyhedron $\polyhedron$ into:

\begin{itemize}
    \item two \emph{principal} convex polyhedra $\polyhedron_1, \polyhedron_2$;
    \item $n$ tetrahedra $\tet^c_{i},\ldots,\tet^c_{i+n}$;
    \item one pyramid $\pir^k_{k-1}$ for each $k\in\{i,\ldots,i+n\}$.
\end{itemize}

\begin{figure}[ht]
    \centering
    \includegraphics[width=.4\linewidth]{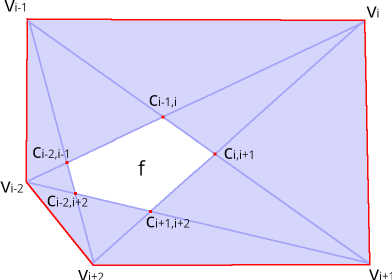}
    \caption{Focus on the face $\face=\polyhedron_1\cap\polyhedron_2$ after chopping all edges: vertices $\pcenter_i^j$ lie on $\plane$, while vertices $\vert_{k}$ are shifted upwards or downwards, according to the sign of $\distance_k$.}
    \label{fig:algorithm:multiple_chop}
\end{figure}

In total, the final set of sub-polyhedra consists of $2n + 2$ polyhedra, of which $2n$ are either tetrahedra or pyramids with square or triangular bases. Let us now focus on the face $\face$ shared by $\polyhedron_1$ and $\polyhedron_2$, which lies on the plane $\plane$ and is shown for this case in Figure~\ref{fig:algorithm:multiple_chop}. Initially bounded by the vertices of $\polygon_0$, this face is progressively modified and ``clipped'' after each chop, giving rise to a new polygon with $n$ vertices $\{\pcenter_i^j\}$. Figure~\ref{fig:algorithm:multiple_chop} shows the face $\face$ after all edges have been chopped, that is, when all distances are nonzero. We note that $\face$ is guaranteed to remain convex because the well-known ear-clipping operation preserves convexity, and it collapses to a single point (the barycenter) when $\polygon_0$ is a quadrilateral.

\subsection{Distances computation}
\label{sec:algorithm:alpha}
This final part is devoted to the discussion about the computation of the distances $\distanceset$ for our algorithm.
As explained at the beginning of this section, each $\distance_i$ is chosen so that the shifted
vertex $\vert_i + \distance_i$ coincides, whenever possible, with one of the endpoints of $\edge_j$ or with its barycenter to reduce the total number of new introduced vertices.
Thus, we first compute the curvilinear coordinate of $\vert_i$ along the edge $\edge_j$ line, measured from the origin of the edge to its endpoint.

In our implementation, if the coordinate belongs to $[0.0, 0.25)$ or to the interval $(0.75, 1.0]$, we compute $\distance_i$ so that $\vert_i$ collapses onto the origin or the endpoint of the edge, respectively.
On the other hand, if the curvilinear coordinate belongs to $[0.25, 0.75]$, we compute $\distance_i$ so that $\vert_i$ coincides with the edge barycenter.
Other choices can be easily implemented to modify the algorithm.

Then, using the computed distance array $\distanceset$, we shift the vertices only if the tetrahedra and the pyramid generated by the shifted vertex have a quality greater than a parameter $\param \in [0,1]$, called \emph{Quality Tolerance}, specified by the user.
The quality of the ``nocht wedges'' can be evaluated using any selected quality metric.
In particular, in our implementation we use the \emph{Scaled Jacobian} indicator, introduced in \cite{KNUPP2003217}.
If $\param = 1$, no quality check is performed; therefore, the \emph{brutal cut} is applied.

\paragraph{Reflection strategy}
\label{sec:reflection}
To avoid the creation of small edges and to further reduce the complexity of the polyhedra, we incorporate a \emph{reflection} strategy into the algorithm.
Specifically, each vertex $\vert_i$ can be shifted not only in one direction, either positive or negative with respect to the normal to the plane $\plane$, but in both directions, thereby reflecting the displacement across the plane, which acts as a ``mirror''.
The algorithm is naturally extended by introducing, for each intersection vertex, the signed distances $\distance^{\pm}_i$, where, in general, $\distance^+_i \neq \distance^-_i$.
The remainder of the algorithm is unchanged, provided that the vertices $\vert_i$ and $\vert_i + \distance_i$ are replaced by $\vert_i + \distance^-_i$ and $\vert_i + \distance^+_i$, respectively.

\section{Results}
\label{sec:results}
\begin{figure}[ht!]
    \centering
    \begin{tabular}{ccccc}
        \includegraphics[height=.12\linewidth]{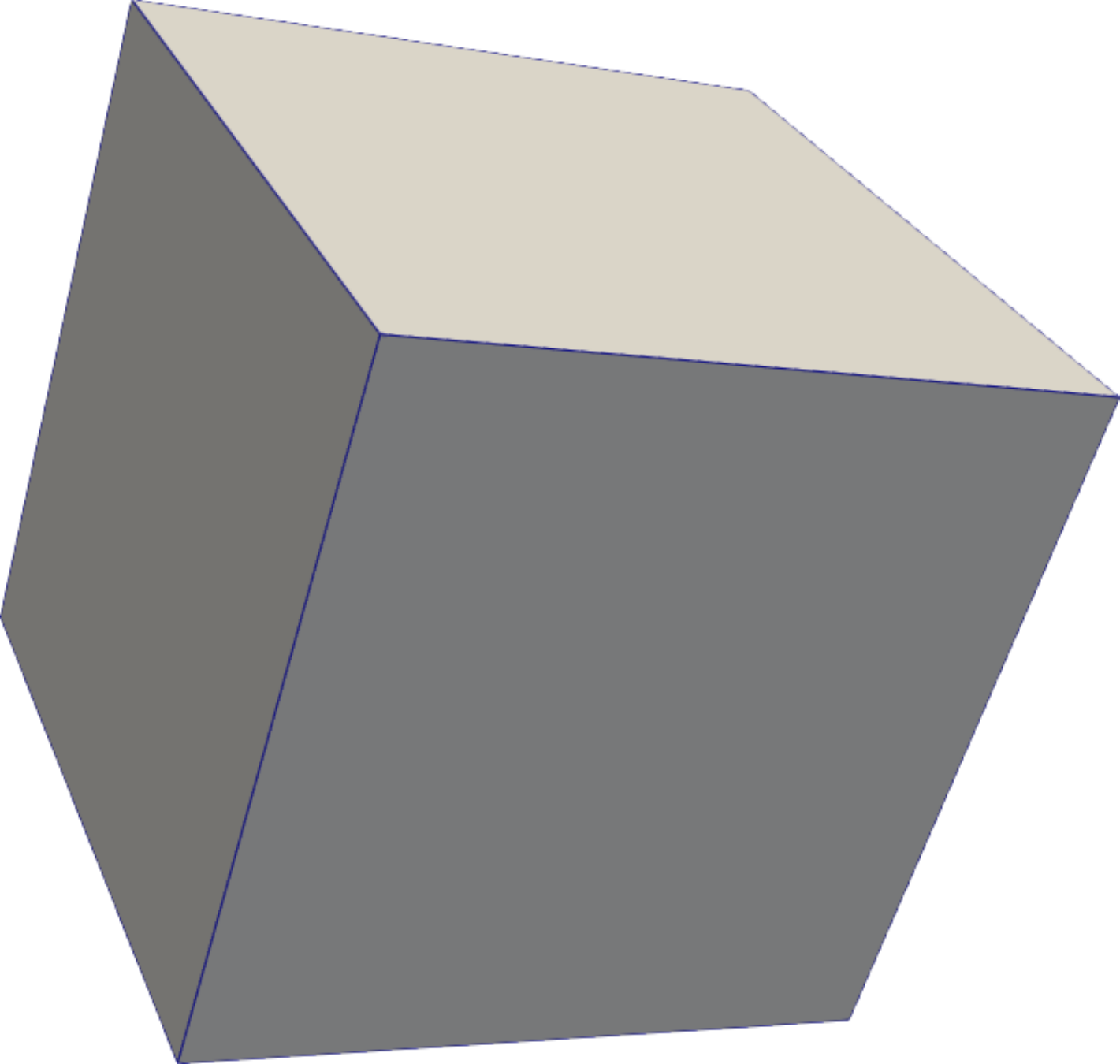} &
        \includegraphics[height=.12\linewidth]{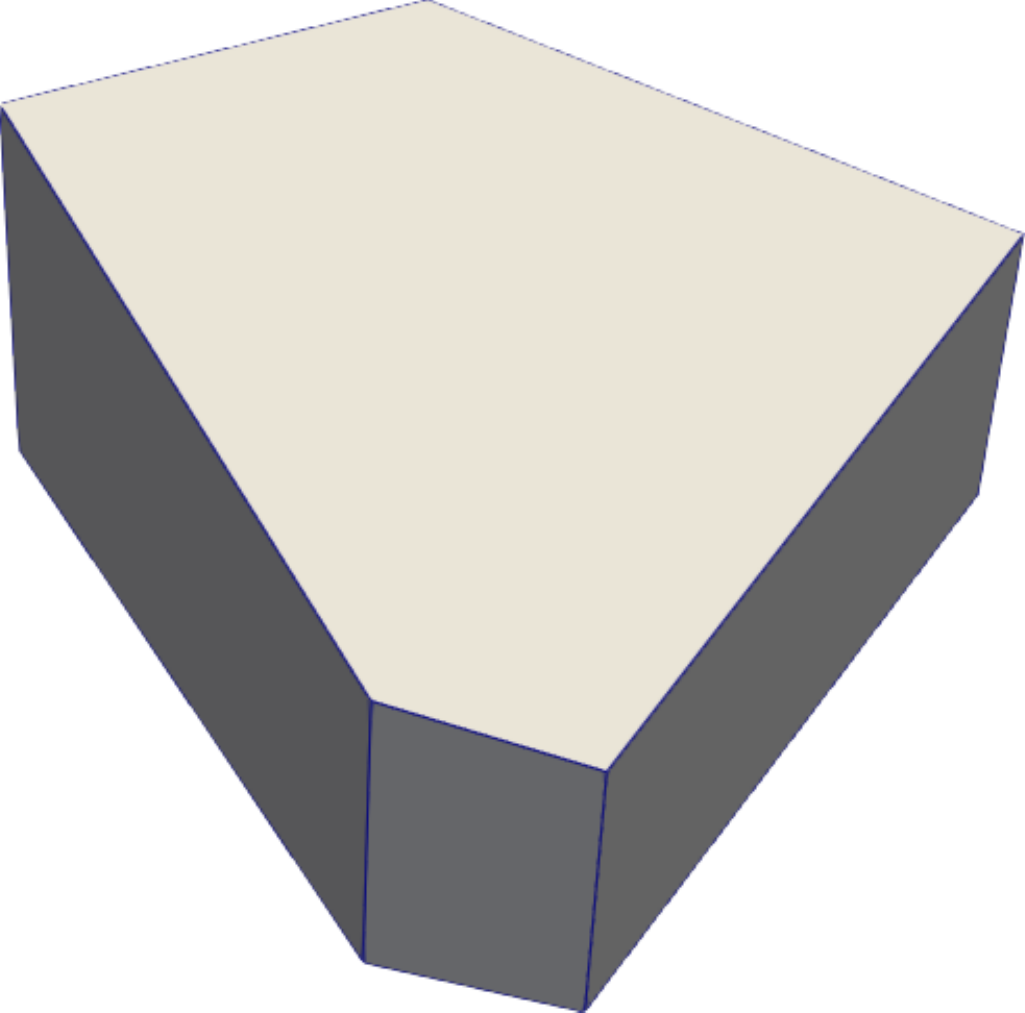} &
        \includegraphics[height=.12\linewidth]{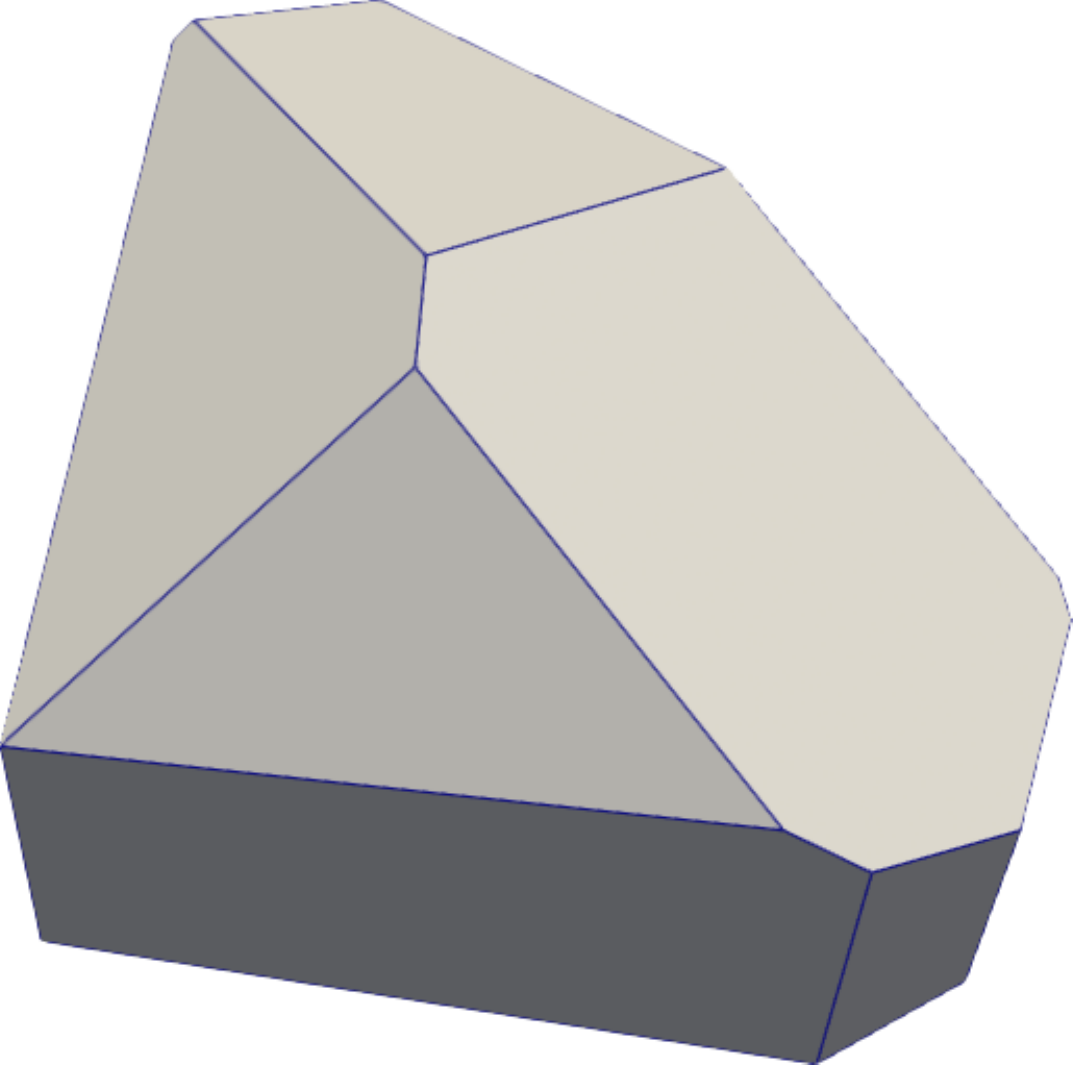} &
        \includegraphics[height=.12\linewidth]{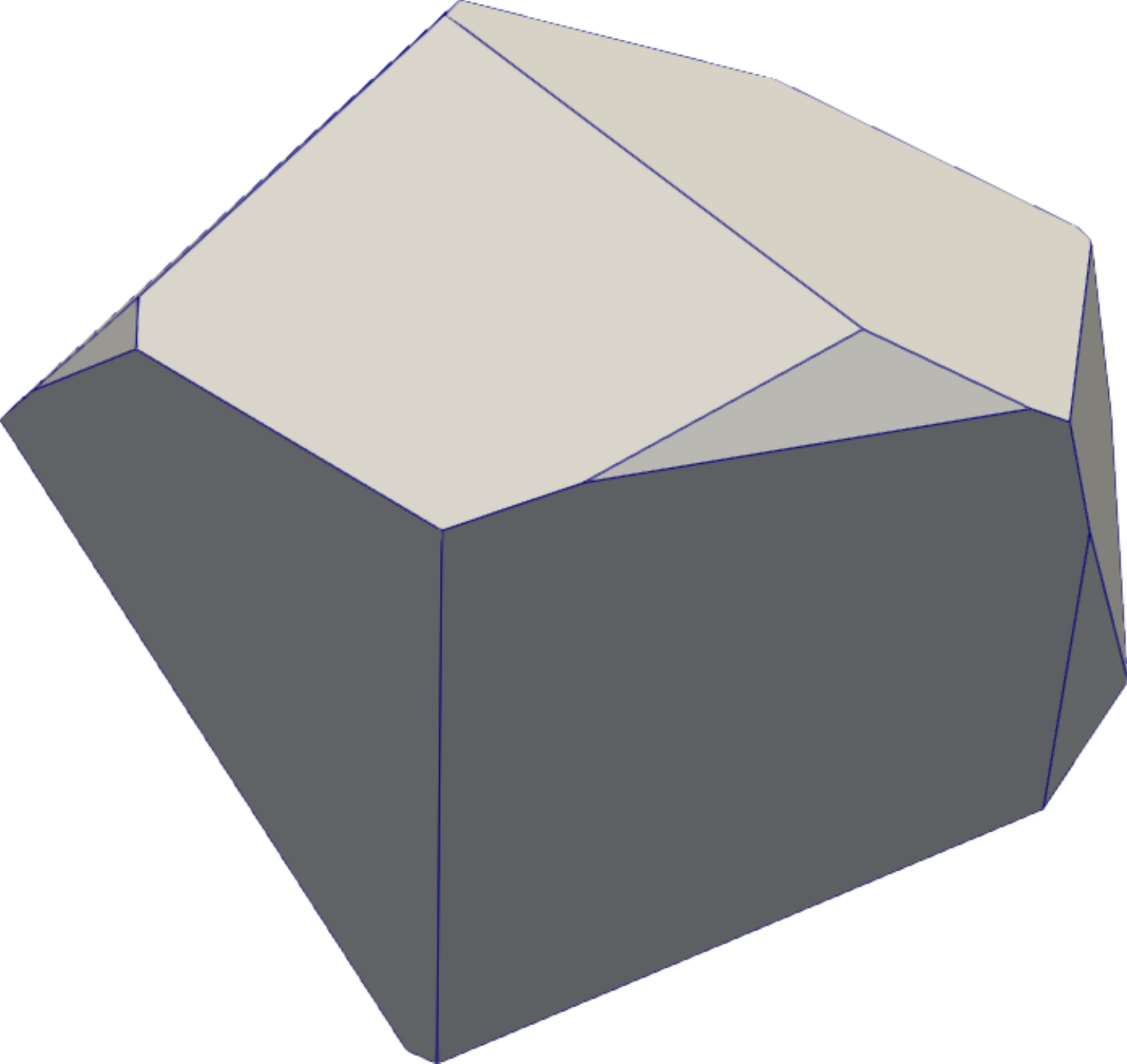} &
        \includegraphics[height=.12\linewidth]{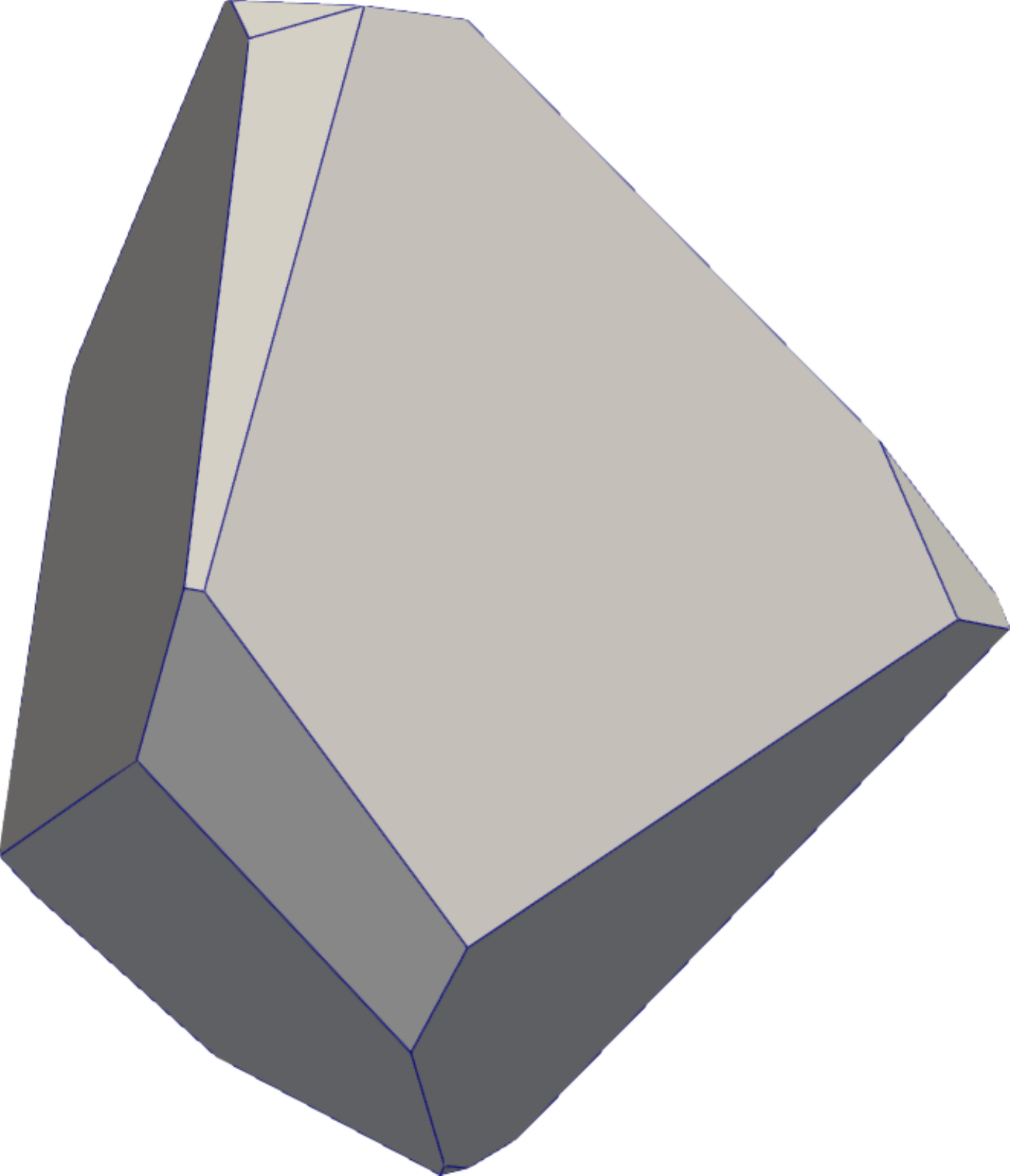} 
    \end{tabular}
    \caption{Examples of polyhedra form our dataset.}
    \label{fig:results:polys}
\end{figure}

We test our algorithm on a sequence of 20 polyhedra with generic polygonal faces and an increasing number of vertices, shown in Figure~\ref{fig:results:polys}.
These polyhedra are generated by using random points as seeds for a constrained Voronoi tessellation of a three-dimensional cubic domain, following the mesh dataset introduced in \cite{SORGENTE2022151}.

We analyze the performance of the algorithm from several perspectives.
The first important aspect is the topological complexity of the resulting decomposition, measured in terms of the number of generated polyhedra, faces, and vertices.
Controlling the complexity of the new polyhedra is crucial to ensure that each splitting step produces configurations that are simpler than the original one. 
When the splitting procedure is applied iteratively, this reduction in complexity promotes a progressive refinement process that ultimately guarantees convergence to tetrahedral elements.

Then, we compare the quality of the generated objects, which can be measured through several indicators, see the survey by Sorgente et. al.~\cite{SBMS23}.
We use an original formulation of the aspect ratio, which measures the ratio between the smallest element (edge of face) and the diameter of the whole polyhedron:
\begin{align}
    \MIE(\polyhedron) := \frac{\min(\min_{\edge\in \polyhedron}\|\edge\|, \min_{\face\in \polyhedron}\sqrt{\|\face\|})}{h_\polyhedron}.
\end{align}

We also adopt the $\VEM$ indicator introduced by Sorgente and colleagues to assess the quality of a mesh with respect to the geometrical requirements of the Virtual Elements Method~\cite{SORGENTE2022151}. It measures the star-shapedness of $\polyhedron$, its aspect ratio, and the number of faces and edges:
\begin{align}
    \VEM(\polyhedron) := \sqrt{\frac{\rho_1(\polyhedron) \rho_2(\polyhedron) + \rho_1(\polyhedron) \rho_3(\polyhedron)}{2}}.
\end{align}
In our context $\rho_1(\polyhedron)$ is always 1, because we are only dealing with convex polyhedra, and
\begin{align*}
    \rho_2(\polyhedron) := \frac{\min(\sqrt[3]{\|\polyhedron\|},\,\min_{\face\in\polyhedron} h_\face)}{h_\polyhedron} +
    \frac{1}{\#\polyhedron} \sum_{\face\in\polyhedron} \frac{\min(\sqrt{\|\face\|}, \ \min_{\edge\in\face} h_\edge)} {h_\face}, \quad
    \rho_3(\polyhedron) := \frac{4}{\#\polyhedron} +
    \frac{1}{\#\polyhedron} \sum_{\face\in\polyhedron}
    \frac{3}{\#\face}.
\end{align*}

\subsection{Chopping one polyhedron}
\label{sec:results:single}
\begin{figure}[ht]
     \centering
     \begin{tabular}{cccc}
         \includegraphics[width=.2\linewidth]{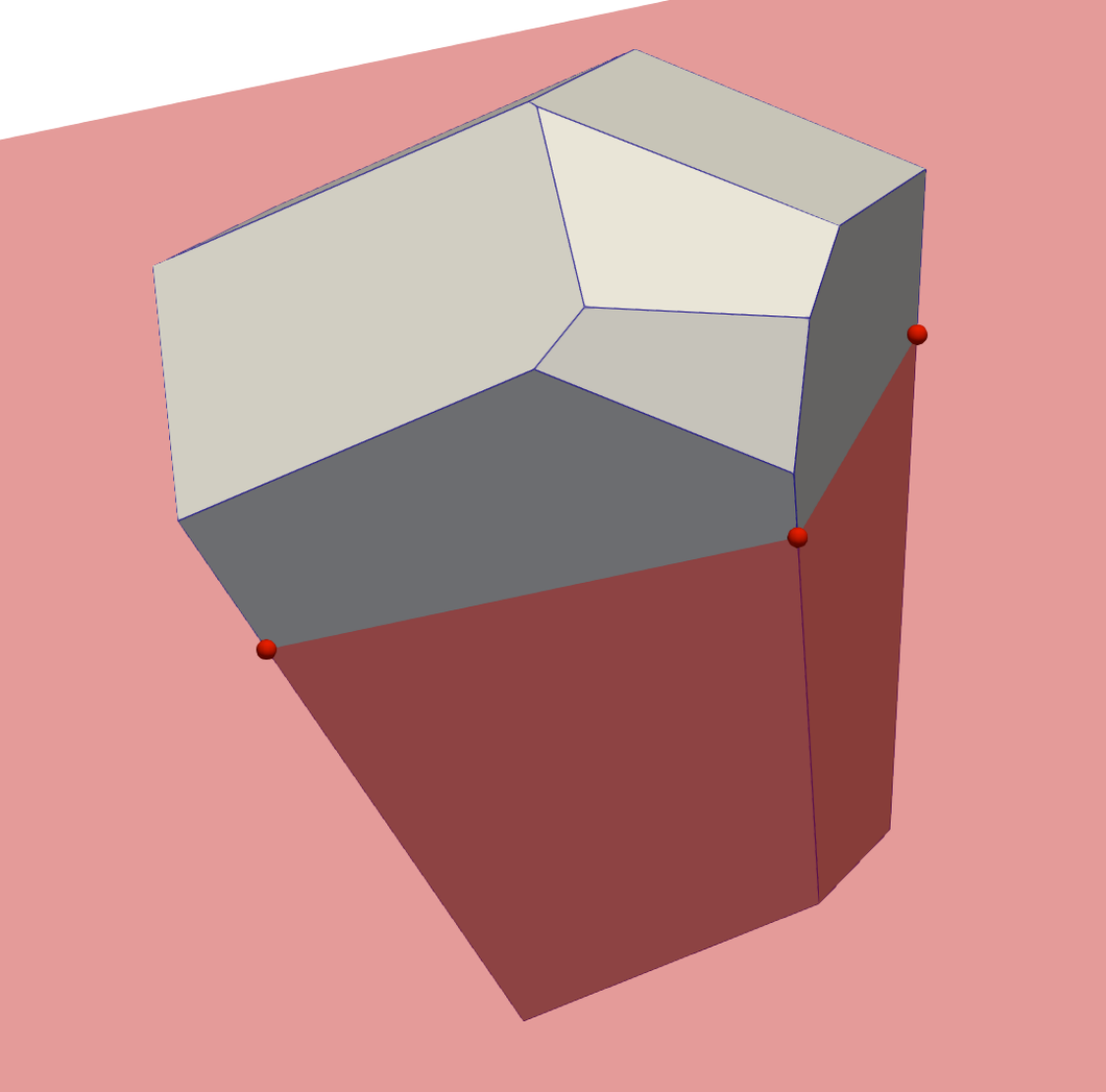} &
         \includegraphics[width=.2\linewidth]{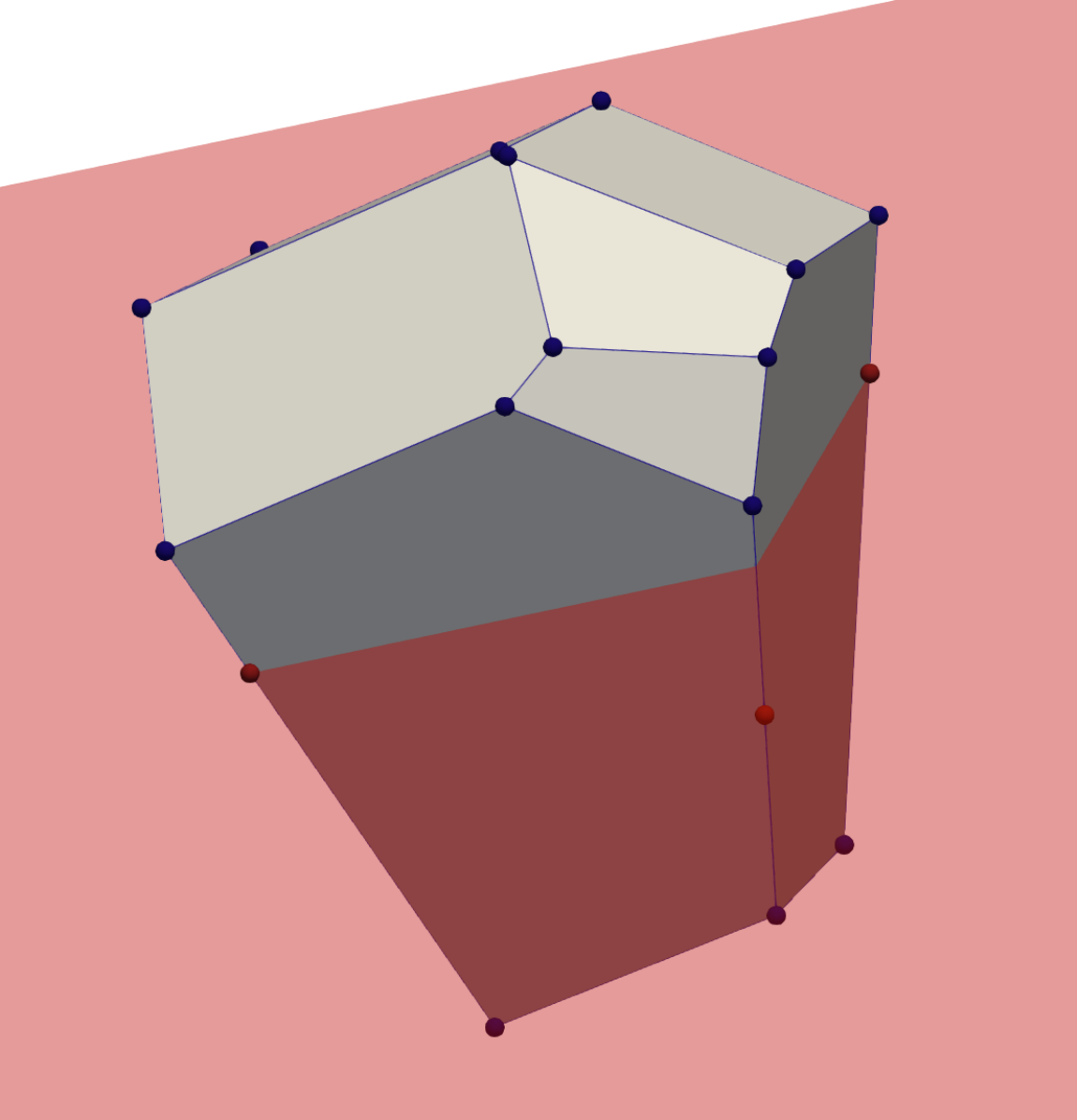} &
         \includegraphics[width=.2\linewidth]{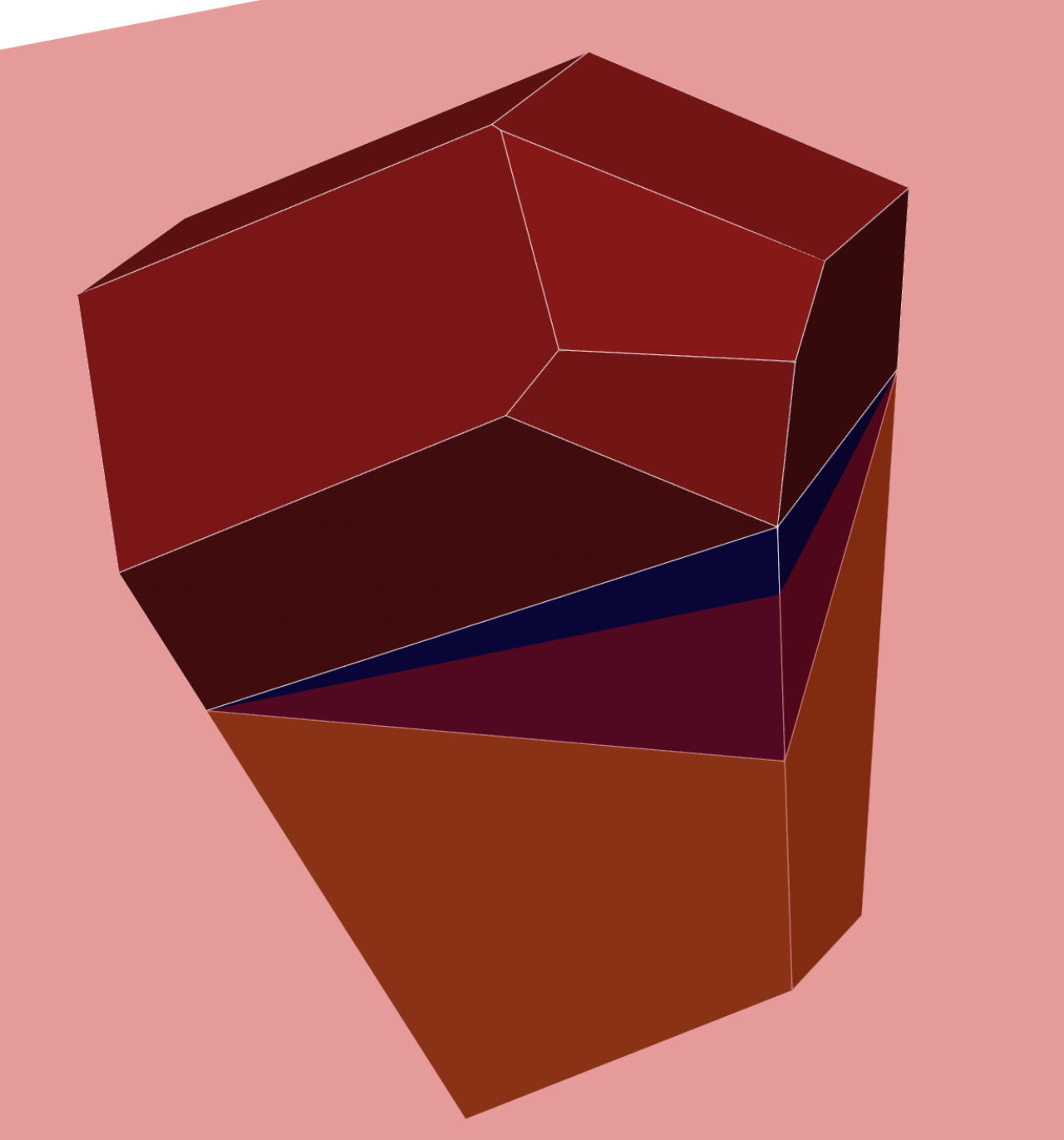} &
         \includegraphics[width=.2\linewidth]{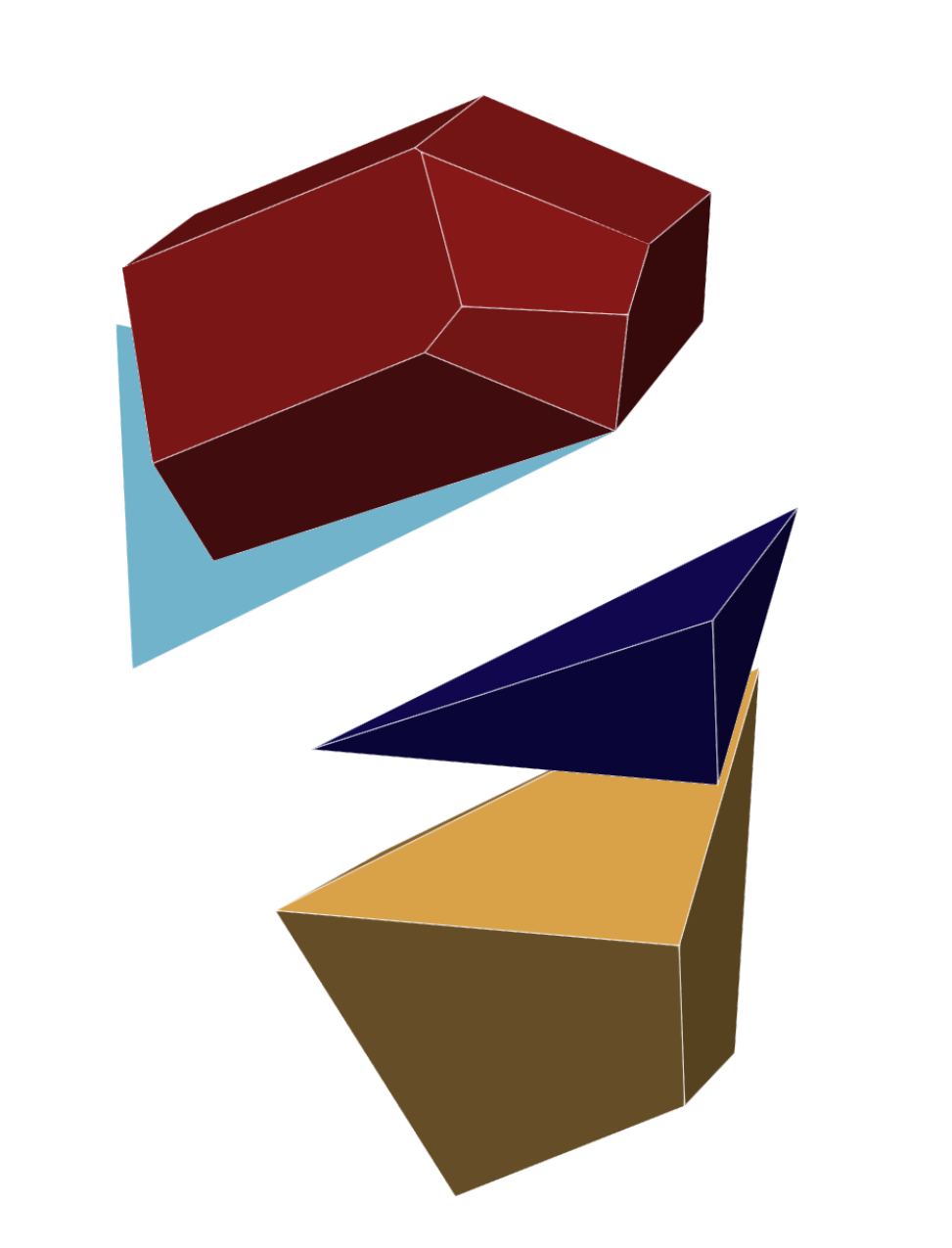} \\
         (a) & (b) & (c) & (d)
     \end{tabular}
     \caption{Chopping of a polyhedron: (a) inertia plane and original intersection points (in red); (b) vertices after chopping: previously existing vertices (in blue) and newly added vertices (in red); (c,d) four sub-polyhedra after chopping.}
     \label{fig:results:single}
\end{figure}

We consider the problem of subdividing a convex polyhedron into smaller convex polyhedra.
In absence of constraints, the PolyChopper algorithm splits the polyhedron along its inertia plane, adjusting the intersection points and using the \textit{reflection} strategy to improve the quality of the resulting objects.
We present this behavior on polyhedron \#3 of our dataset, chopped with parameter $\param=0.4$.
Figure~\ref{fig:results:single} illustrates how some of the original intersection points with the inertia plane are moved towards the middle of the respective edges, to prevent the generation of small edges and faces.
Four new objects are generated: two \textit{main} sub-polyhedra (in red and orange) and two \textit{wedges} sub-polyhedra (in blue and light blue).

\begin{figure}[ht]
     \centering
     \includegraphics[width=.75\linewidth]{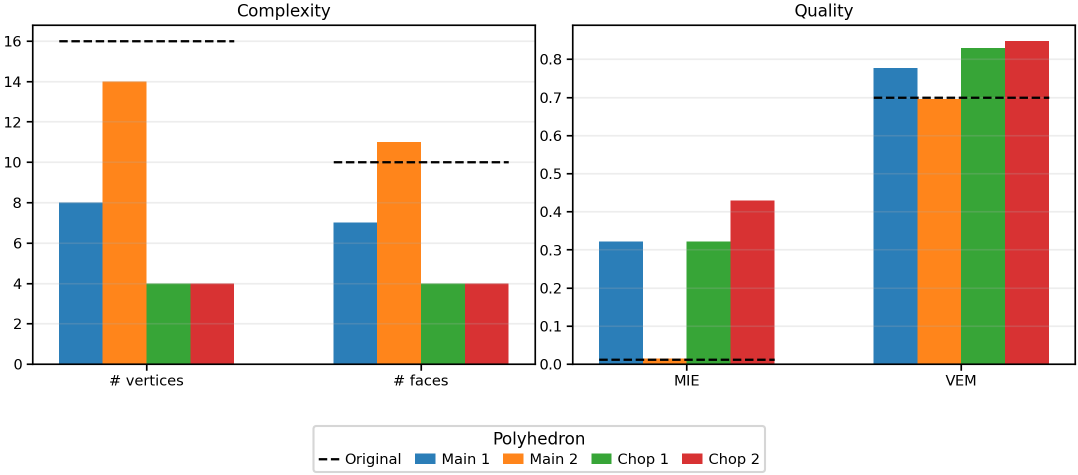}
     \caption{Complexity and quality analysis for chopping of polyhedron \#3 with $\param=0.4$.}
     \label{fig:results:single-stats}
\end{figure}

Figure~\ref{fig:results:single-stats} presents the complexity and quality analysis.
The sub-polyhedron \textit{Main 2} (red in Figure~\ref{fig:results:single}) is comparable to the original polyhedron, having less vertices but slightly more faces, and also similar quality scores.
However, all other sub-polys are significantly simpler, and have higher quality for both indicators.
In conclusion, our algorithm managed to subdivide this polyhedron into three simpler and higher quality sub-elements, plus one element of comparable complexity and quality.

\begin{figure}[ht]
     \centering
     \includegraphics[width=.8\linewidth]{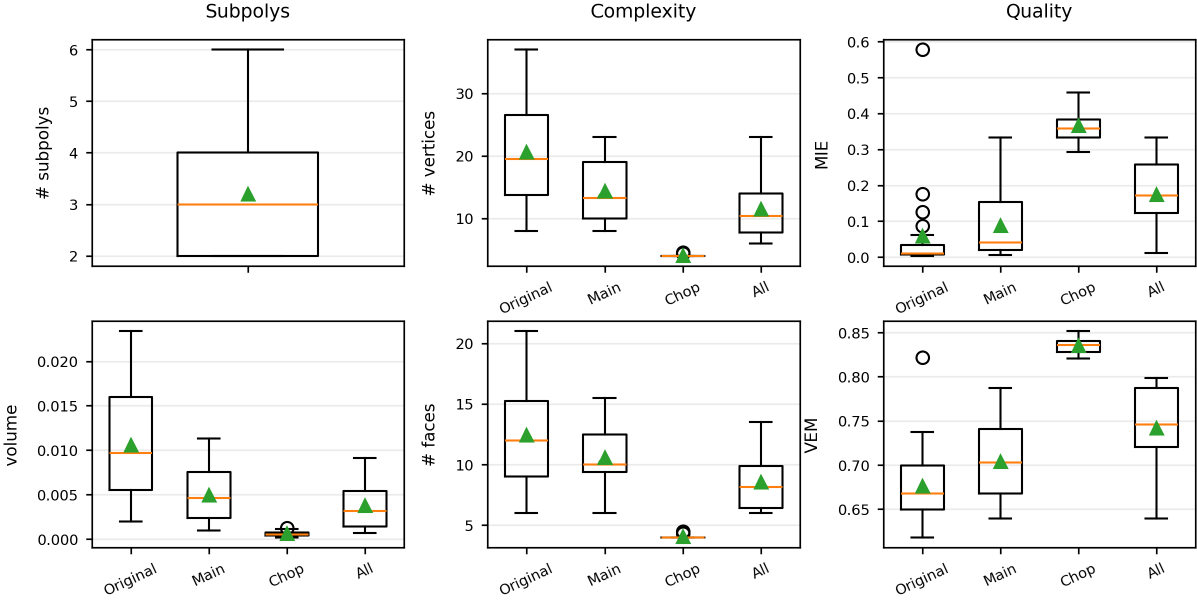}
     \caption{Complexity and quality statistics for the whole dataset with $\param=0.4$. 
     Boxes span the interquartile range ($Q_1$--$Q_3$), the horizontal line marks the median, the triangle marks the mean, whiskers indicate the non-outlier range, and isolated points are outliers.}
     \label{fig:results:generic_dataset}
\end{figure}

To broaden our analysis, we perform the same chopping on all elements from the dataset, using $\param=0.4$; results are collected in Figure~\ref{fig:results:generic_dataset}.
\emph{Original} refers to the input polyhedron, while \emph{Main}, \emph{Chop}, and \emph{All} refer to the mean values among the two main output sub-polyhedra, the small chopped sub-polyhedra, and all generated sub-polyhedra, respectively.
We observe how the median number of sub-polys generated by the algorithm is 3, meaning that most times we have only one extra chop besides the two main polyhedra.
The decomposition systematically reduces the combinatorial complexity of the generated sub-polyhedra with respect to the original polyhedra. This reduction is visible in the number of vertices and faces. The two main output polyhedra are still the most complex generated objects, but their median complexity is lower than that of the original polyhedra. The chopped elements are instead much simpler, with very small numbers of vertices and faces, as expected from their restricted shape classes.
In terms of quality, the main sub-polyhedra remain close to the original distribution, improving the median values of both MIE and VEM.
This suggests that the decomposition reduces complexity while preserving, and in some cases mildly improving, the quality of the main generated components. 

\subsection{Role of the parameter $\param$}
\label{sec:results:parameter}
\begin{figure}[ht]
     \centering
     \includegraphics[width=.8\linewidth]{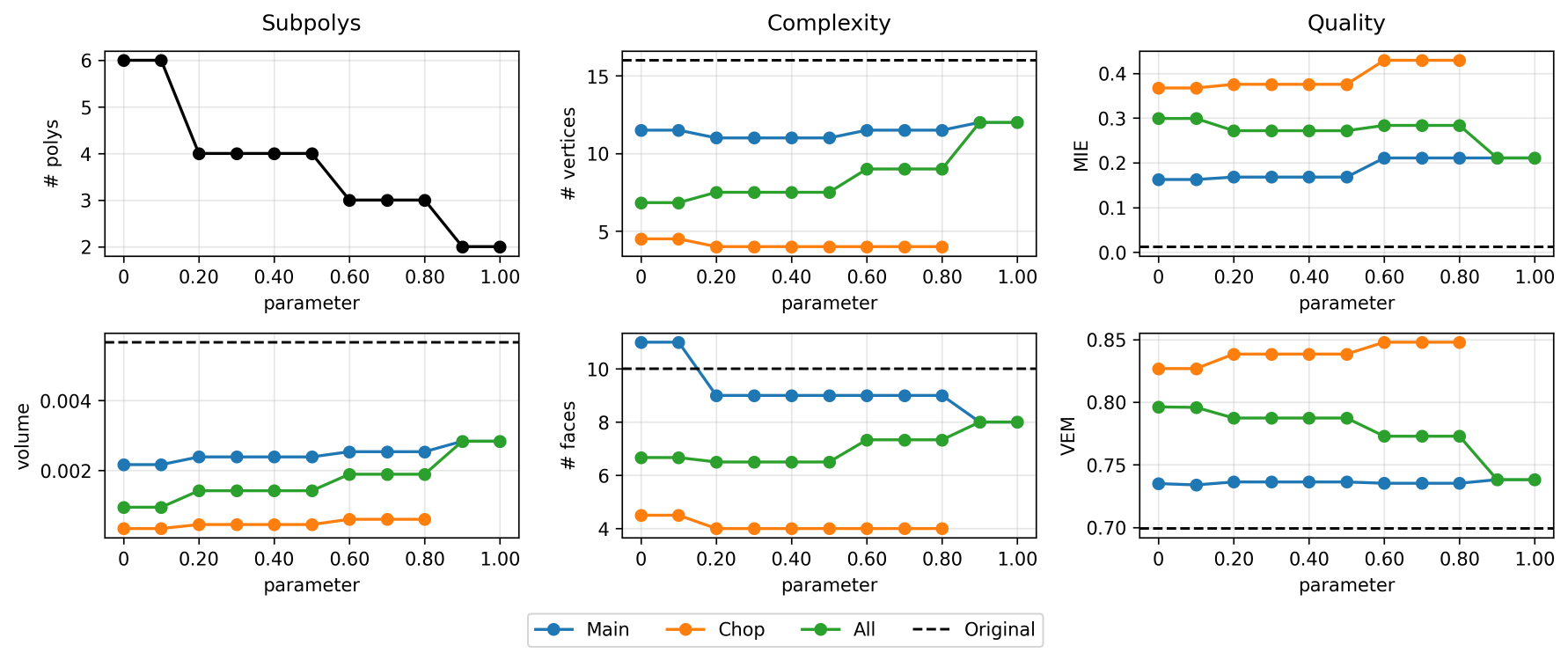}
     \caption{Complexity and quality analysis for chopping of polyhedra \#3 with varying parameter $\param\in[0,1]$.}
     \label{fig:results:angle-stats}
\end{figure}

The \emph{Quality Tolerance} parameter $\param$, introduced in Section~\ref{sec:algorithm:alpha}, has a key role in the algorithm, as presented in Figure~\ref{fig:results:angle-stats}.
We recall thatm, if $\param=1$, the algorithm performs the \textit{brutal cut}: the polyhedron is divided by the plane without moving any point and without generating wedges, i.e., it produces only two main sub-polys.
As $\param\to 0$, more and more intersection points are allowed to be moved to increase the quality of the two main sub-polys, generating chop polyhedra.

The mean complexity of the two main polyhedra remains approximately constant, as the number of vertices decreases but the number of faces increases, except for $\param\le0.1$.
For $\param\in[0.6,0.8]$ some isolated tetrahedra appear from single chops, see Figure~\ref{fig:results:single}(d), while for $\param\le0.5$ we have pyramids generated by consecutive chops.
All these elements are way less complex than the original polyhedron, contributing to a decrease of the average complexity.
The quality of the main sub-polyhedra also remains constant, but the higher quality of the chop sub-polys increases the mean quality score,  with tetrahedral chops having higher quality with respect to pyramidal ones.

Overall, cutting the polyhedron along its inertia plane produces simpler and better sub-polyhedra by itself, but lower values of $\param$ allow to introduce some chop sub-polys that decrease the average complexity and increase the average quality of the resulting decomposition.

\subsection{Iterative chopping for mesh generation}
\label{sec:mesh}
\begin{figure}[ht]
     \centering
     \begin{tabular}{ccc}
         \includegraphics[width=.22\linewidth]{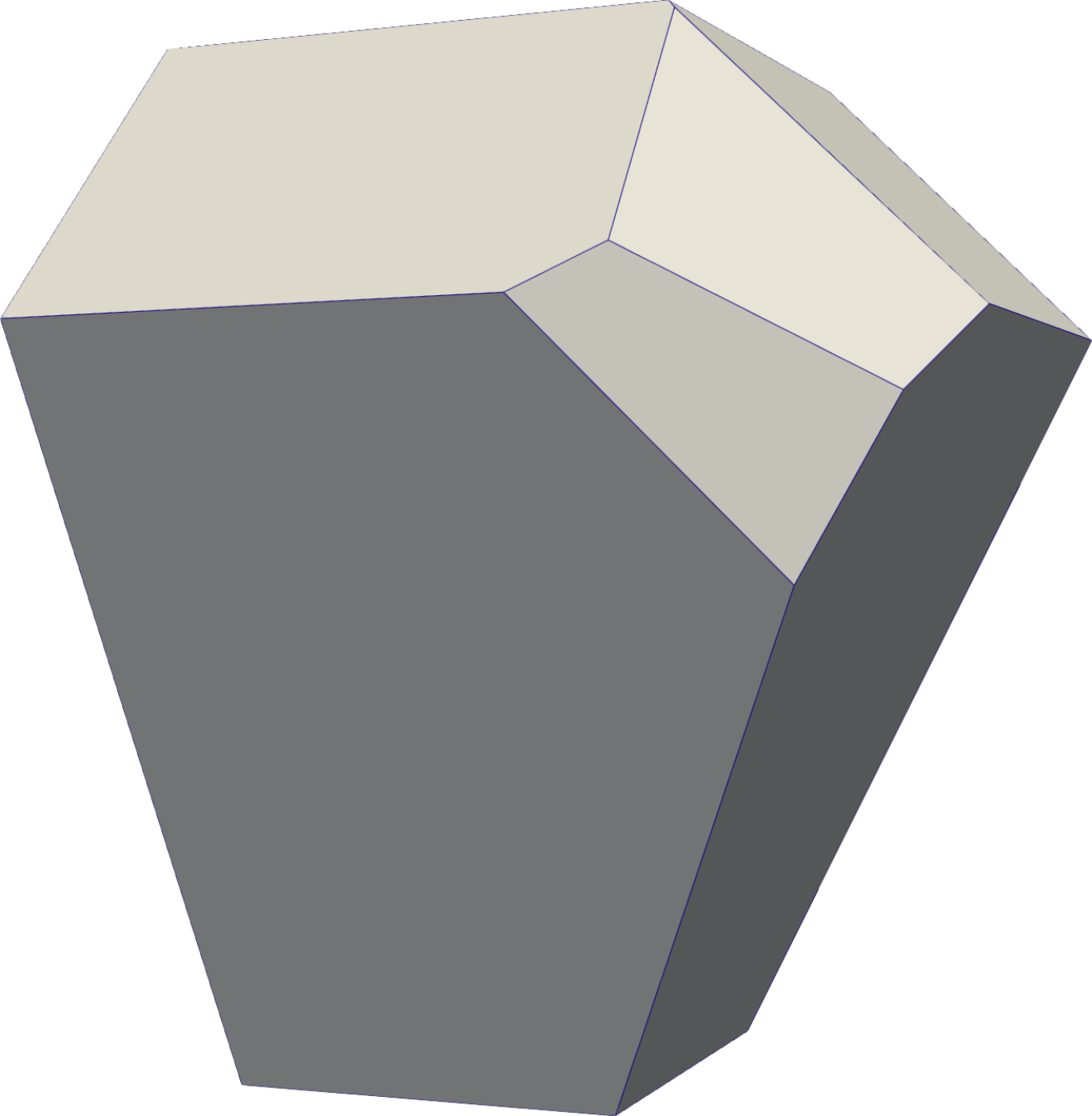} &
         \includegraphics[width=.22\linewidth]{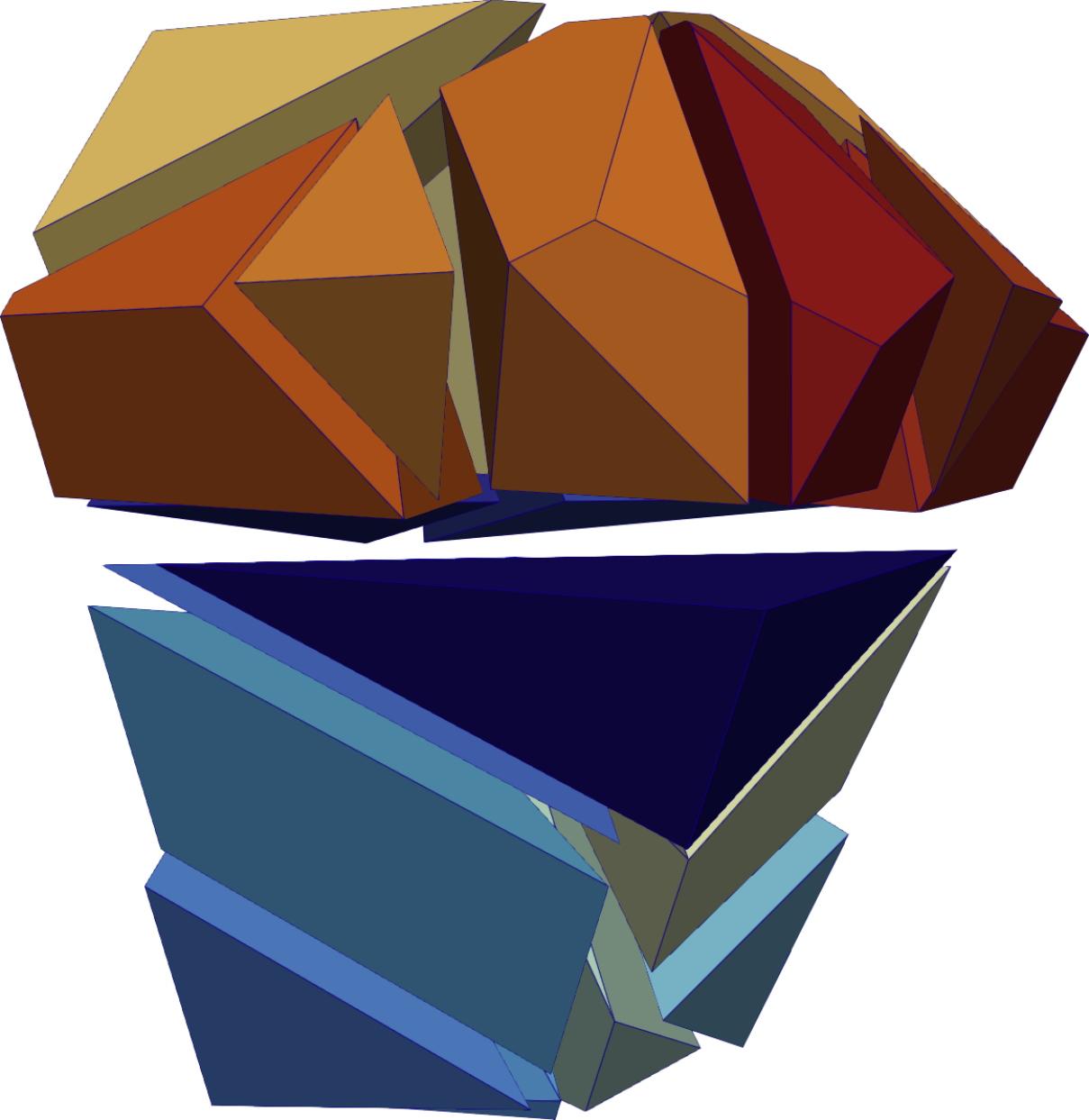} &
         \includegraphics[width=.22\linewidth]{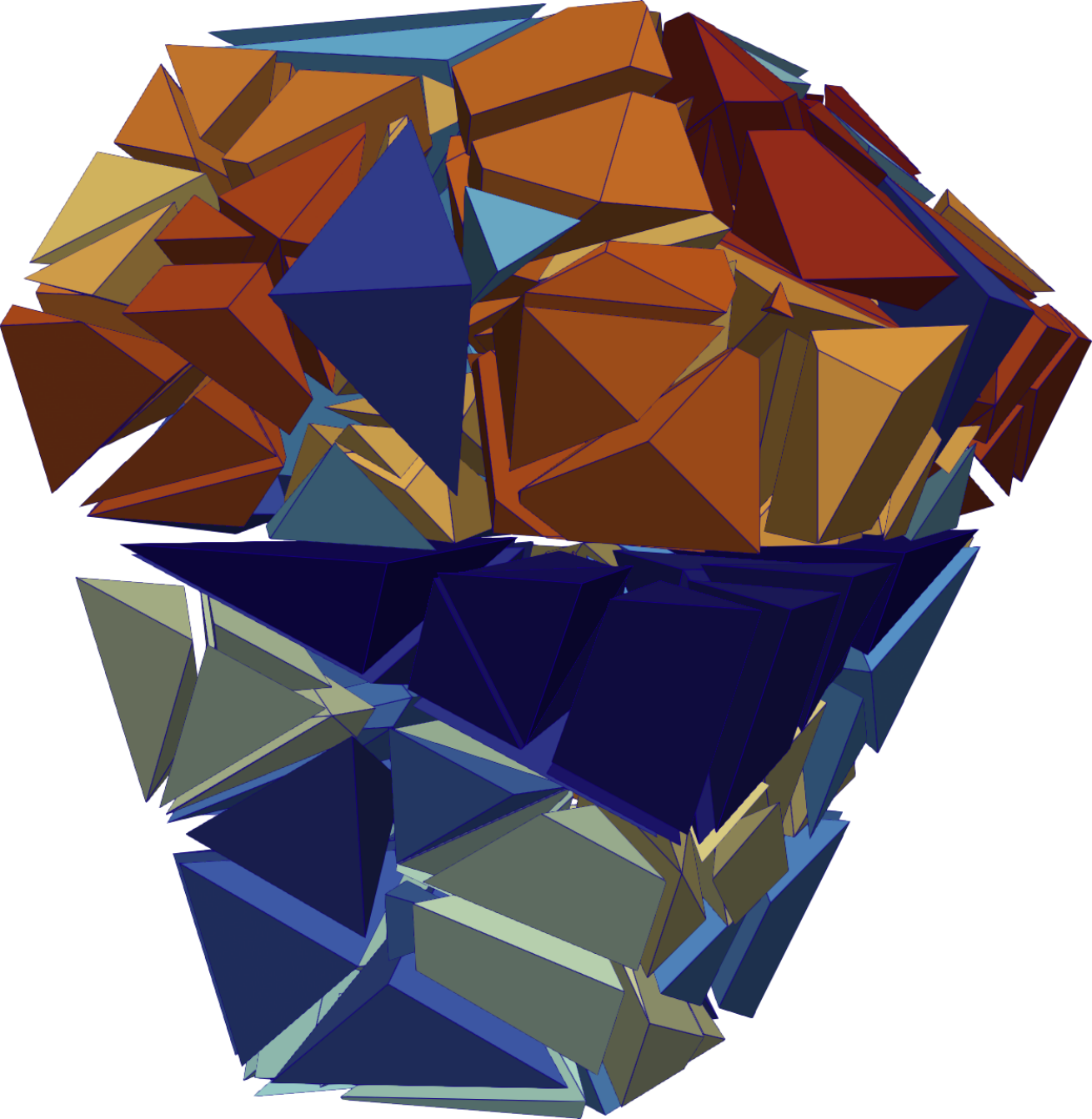} \\
         (a) & (b) & (c)
     \end{tabular}
     \caption{Iterative chopping of polyhedron \#3 (a) with $\Vparam=0.1$ (b) and $\Vparam=0.01$ (c).}
     \label{fig:results:refinement_mesh}
\end{figure}

One interesting application of our algorithm is within meshing algorithms.
Indeed, chopping can be performed iteratively to refine the cells of a polyhedral mesh, or to generate a mesh of a convex domain.
To this aim, we introduce an apposite parameter $\Vparam$ which regulates the relative volume reduction of the cells: $\Vparam=0.1$ means that the volume of each sub-polyhedron will be less than 10\% of the volume of the original polyhedron.
In Figure~\ref{fig:results:refinement_mesh}, we present the iterative chopping of polyhedron \#3 with different $\Vparam$ values.
This generates a sequence of encapsulated subdivisions of the original polyhedron, which can also be employed for multi-resolution applications.

\begin{figure}[ht]
     \centering
     \begin{tabular}{ccc}
         \includegraphics[width=.3\linewidth]{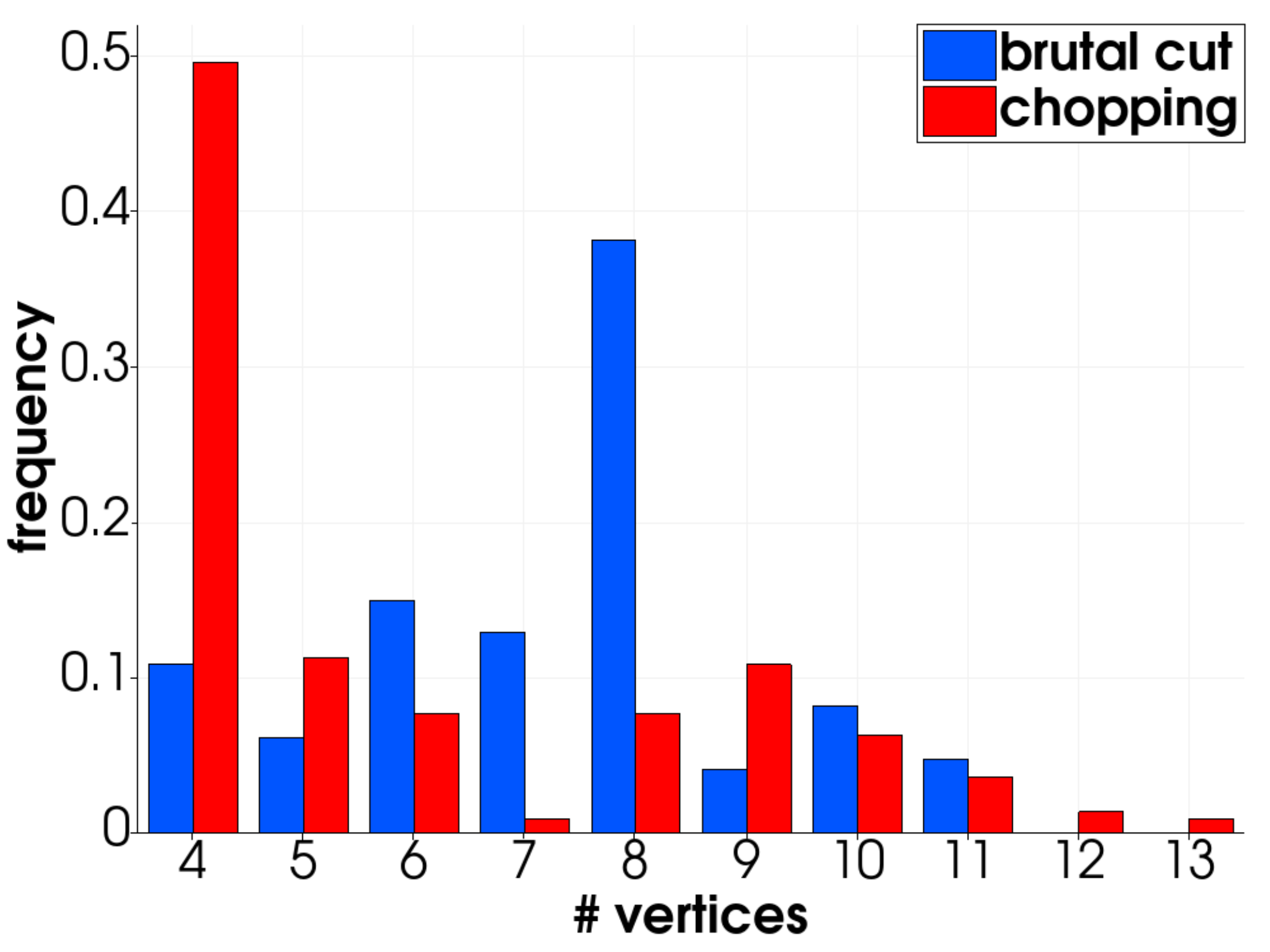} &
         \includegraphics[width=.3\linewidth]{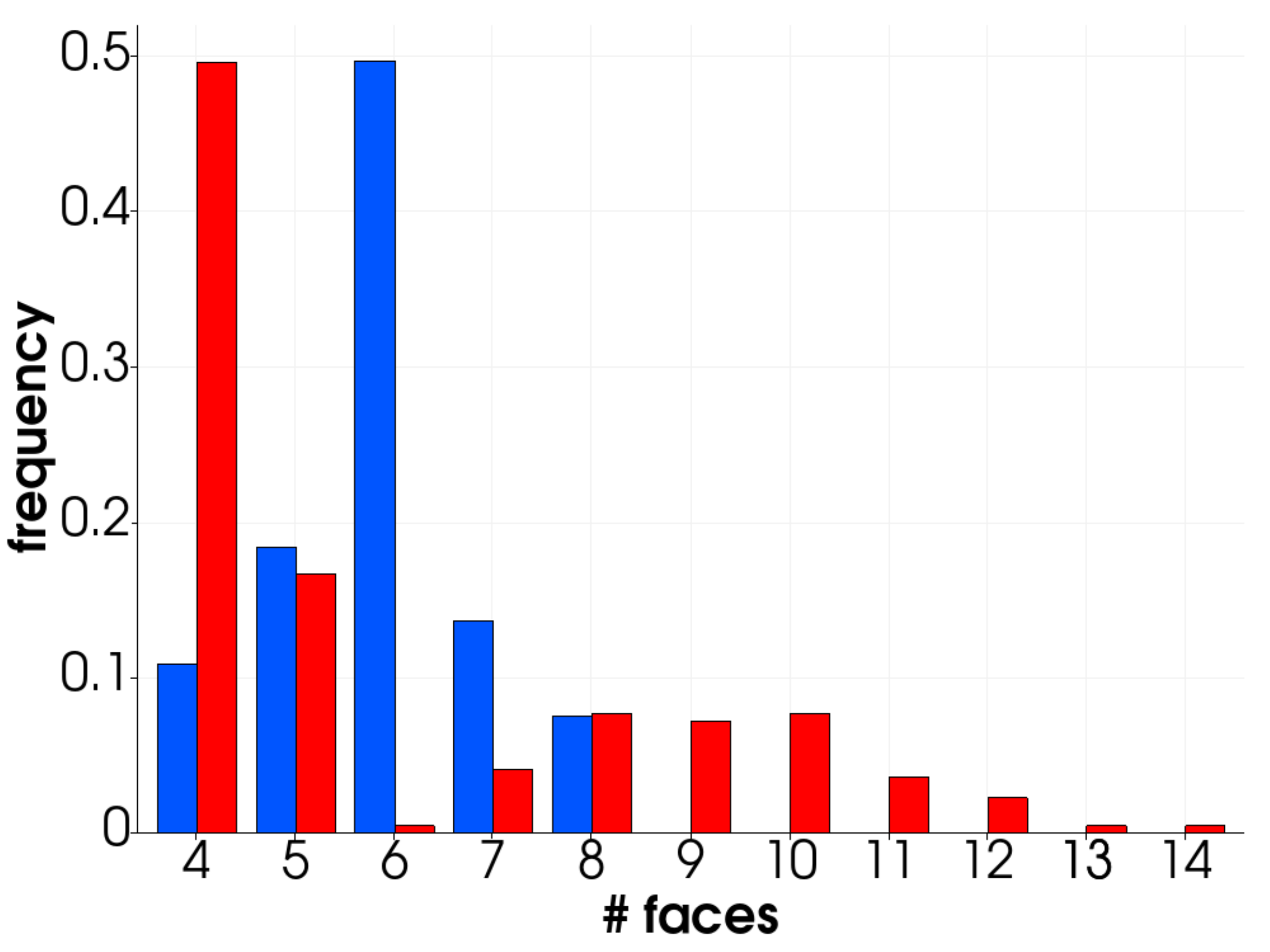} &
         \includegraphics[width=.3\linewidth]{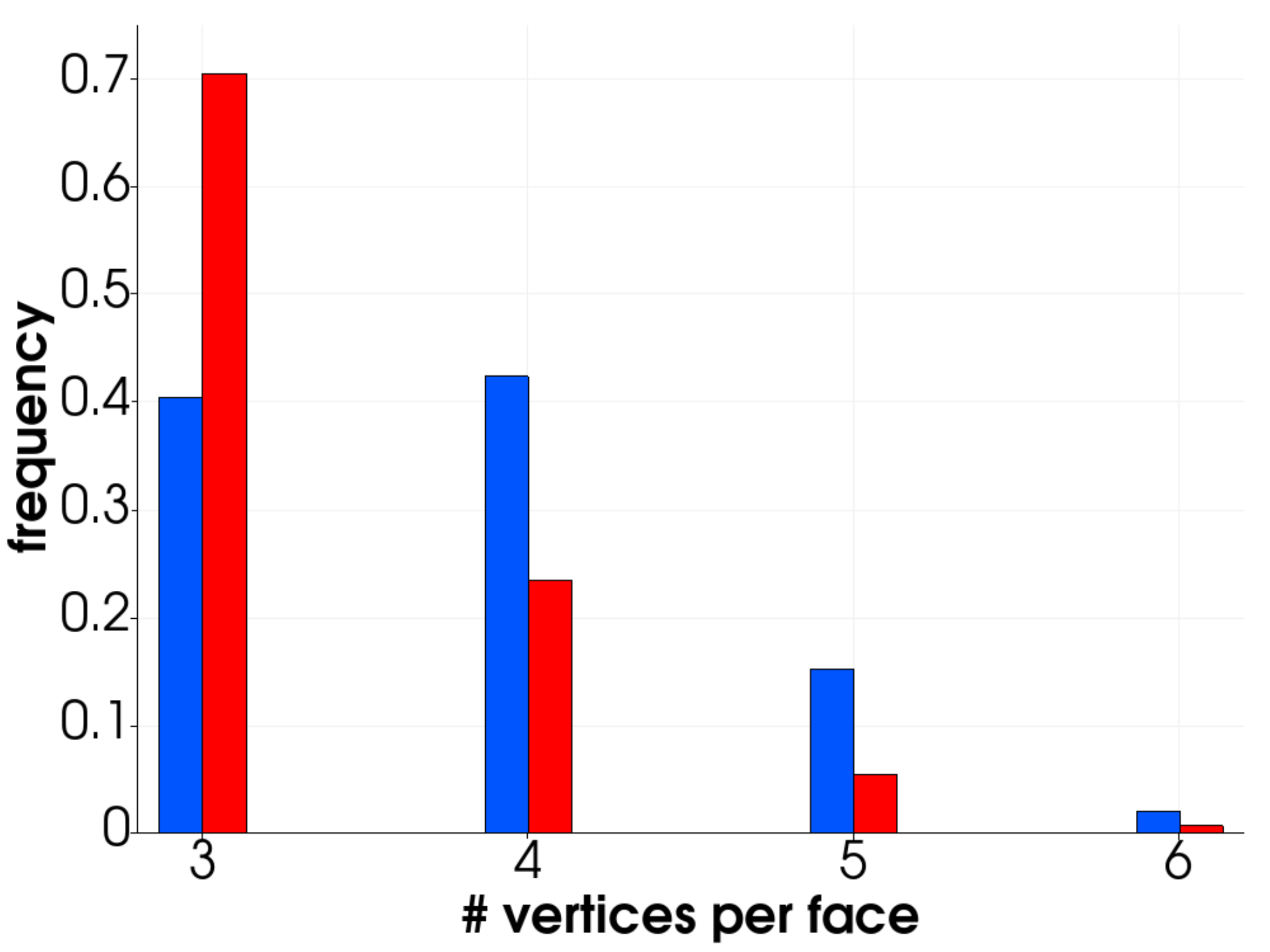} 
     \end{tabular}
     \begin{tabular}{cc}
         \includegraphics[width=.3\linewidth]{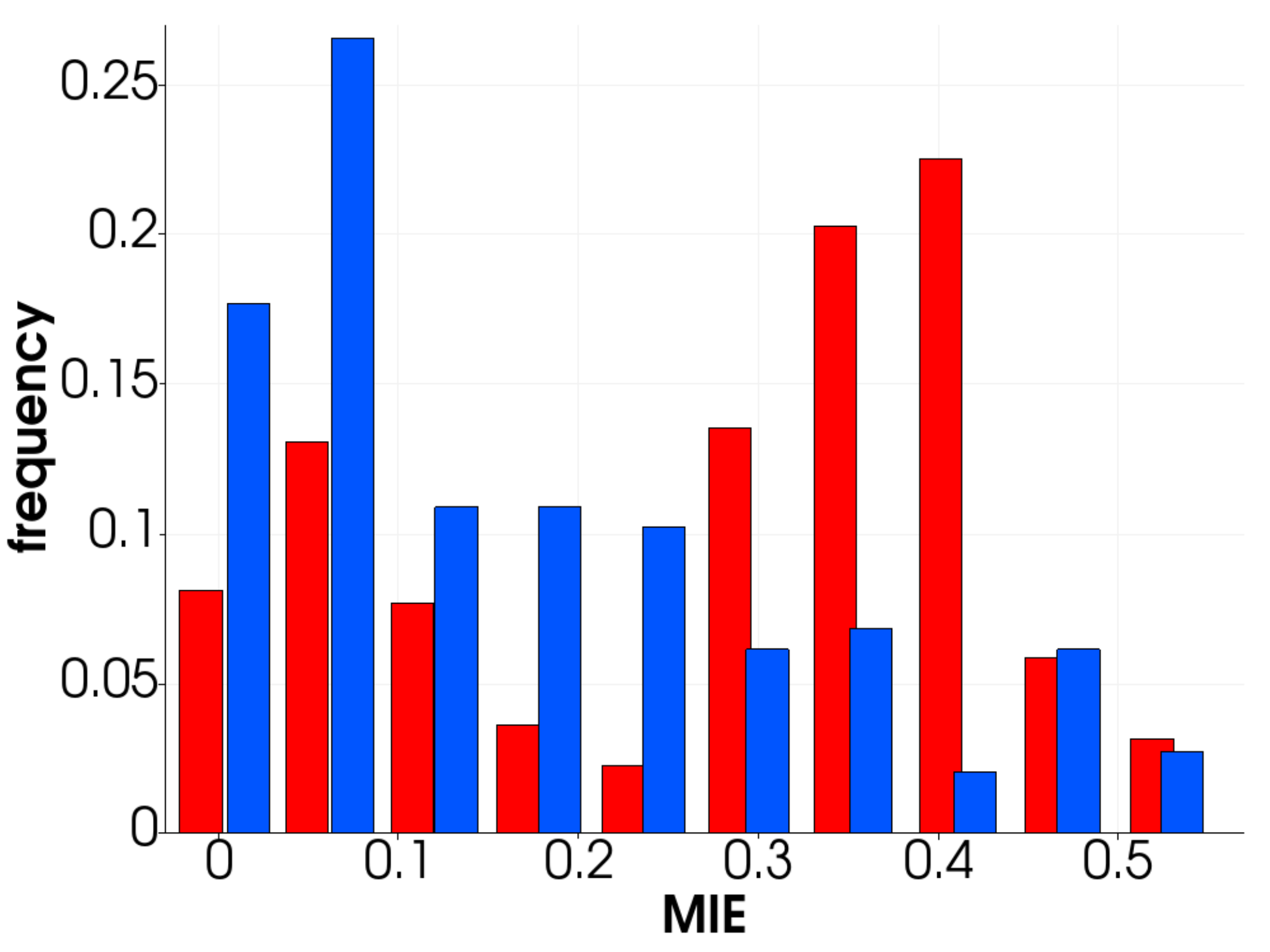} &
         \includegraphics[width=.3\linewidth]{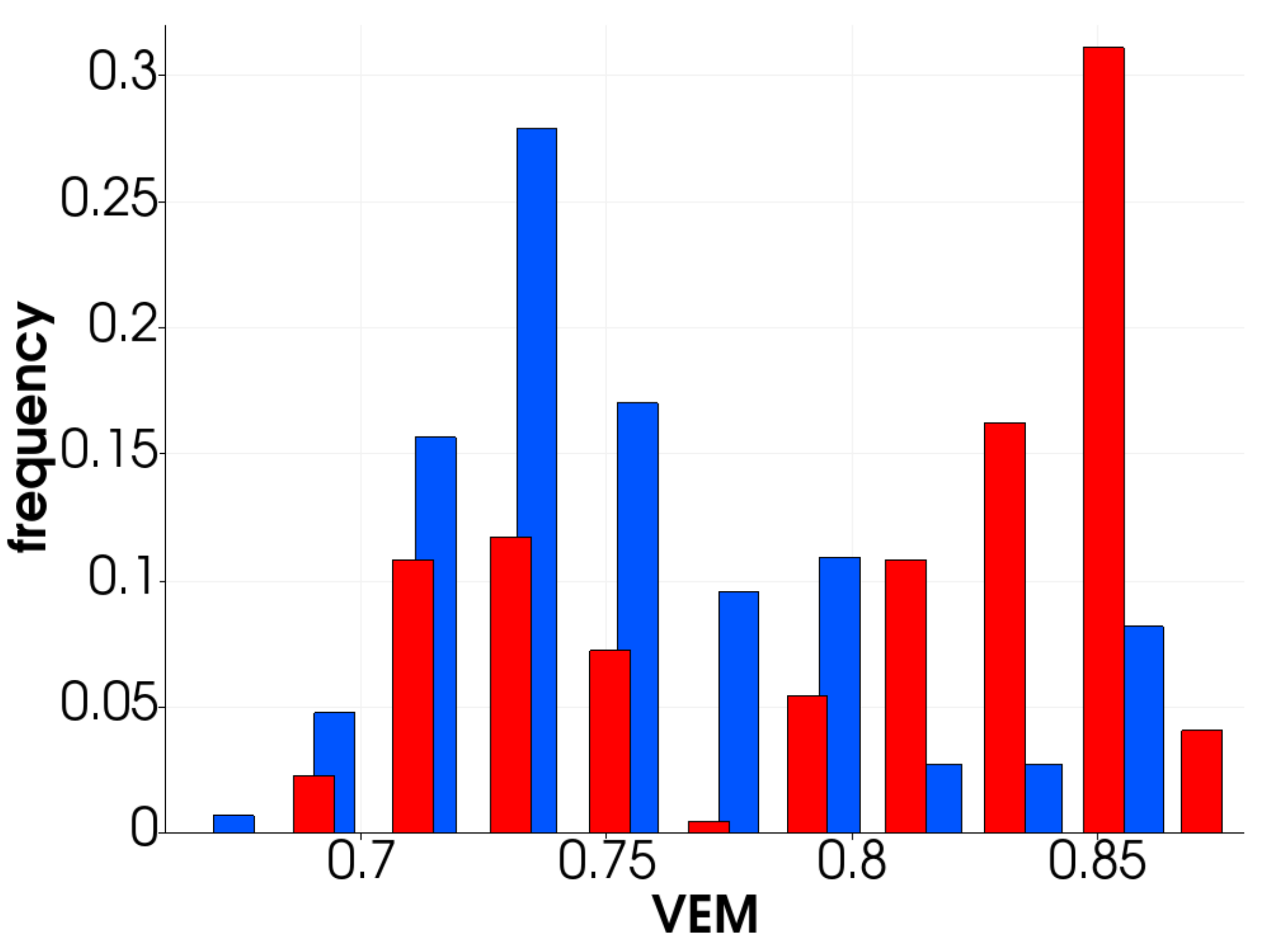} 
     \end{tabular}
     \caption{Statistics about the iterative subdivision of polyhedron \#3 with $\Vparam=0.01$, using the brutal cut and the PolyChopper algorithm.}
     \label{fig:results:refinement_stats}
\end{figure}

Figure~\ref{fig:results:refinement_stats} compares the iterative subdivision of a polyhedron using the same value $\Vparam=0.01$ and two different techniques, namely, the PolyChopper algorithm with $\param=0.4$ and the brutal cut approach.
The latter is obtained by setting $\param=1$, and it is used as a comparison to understand the benefits of the introduction of chops in the subdivision.
A major difference is that chopping generates tet-dominant meshes while the brutal cut does not: this is appreciable by observing the number of sub-polys with 4 vertices (equivalently, 4 faces).
Besides a 50\% of tetrahedral elements, chopping generates a 70\% of triangular faces, that is, even on non-tetrahedral elements most faces are triangular.
Elements generated by brutal cuts are, on average, more complex, and with lower quality.
Indeed, for both $\MIE$ and $\VEM$, the element distributions obtained with the chopping algorithm are shifted toward higher values, and the distribution peaks occur at significantly larger values than those associated with the brutal cut method.
This indicates the superiority of the subdivision generated by our algorithm, and we note that, after the first few cuts, the generated sub-polyhedra become largely independent of the original polyhedron. 
Therefore, this result can be generalized to any other original polyhedron, as confirmed by our experiments.

\section{Conclusions}
\label{sec:conclusions}
We presented \textit{PolyChopper}, an algorithm for decomposing a convex polyhedron into a finite set of convex sub-polyhedra. 
The decomposition is achieved by cutting the polyhedron with a plane, either prescribed as an external constraint or automatically determined from the inertia tensor. The intersection vertices are then adjusted to control both the complexity and the quality of the resulting elements, while additional cuts are introduced according to a quality parameter. 
Experimental results show that the proposed approach is robust and consistently generates sub-polyhedra that are simpler and of higher quality than the original element. Furthermore, the recursive decomposition produces a hierarchy of polyhedral subdivisions with progressively smaller element sizes, while preserving convexity and generating polyhedra that are predominantly tetrahedral.

To the best of our knowledge, no other existing algorithm produces polyhedral tessellations with the same combination of properties. Future work will include a more extensive comparison with other polyhedral meshing techniques, such as Voronoi-based tessellations, as well as numerical studies using polytopal discretization methods, including the Virtual Element Method. 
We expect the tetrahedral-dominant meshes generated by \textit{PolyChopper} to be particularly effective for domains with complex boundaries or internal constraints, where their flexibility enables the accurate representation of small-scale geometric features without a significant increase in mesh complexity (i.e., the number of elements, vertices, and faces). Another promising direction is the extension of the method to support non-planar constraint surfaces.

\section*{Acknowledgments}
T. Sorgente and F. Vicini are affiliated to the Italian Gruppo Nazionale Calcolo Scientifico - Istituto Nazionale di Alta Matematica (GNCS-INdAM).
T. Sorgente acknowledges the financial support provided by INdAM-GNCS Project ``Metodi numerici politopali stabilization-free e neural-based per problemi accoppiati e non lineari'' (CUP\_E53C25002010001).
F. Vicini acknowledges the financial support provided by INdAM-GNCS Project ``Metodi numerici integrati per la simulazione e la prevenzione del dissesto idrogeologico'' (CUP\_E53C25002010001).

\printbibliography

@article{SBMS23,
author = {Sorgente, T. and Biasotti, S. and Manzini, G. and Spagnuolo, M.},
title = {A Survey of Indicators for Mesh Quality Assessment},
journal = {Computer Graphics Forum},
volume = {42},
number = {2},
pages = {461-483},
doi = {https://doi.org/10.1111/cgf.14779},
year = {2023}
}

@article{10.1007/s10092-023-00517-5,
author = {Sorgente, Tommaso and Vicini, Fabio and Berrone, Stefano and Biasotti, Silvia and Manzini, Gianmarco and Spagnuolo, Michela},
title = {Mesh quality agglomeration algorithm for the virtual element method applied to discrete fracture networks},
year = {2023},
issue_date = {Jun 2023},
publisher = {Springer-Verlag},
address = {Berlin, Heidelberg},
volume = {60},
number = {2},
issn = {0008-0624},
url = {https://doi.org/10.1007/s10092-023-00517-5},
doi = {10.1007/s10092-023-00517-5},
journal = {Calcolo},
month = apr,
numpages = {27}
}

@article{berrone2023optimization,
  title={An optimization based 3D-1D coupling strategy for tissue perfusion and chemical transport during tumor-induced angiogenesis},
  author={Berrone, Stefano and Giverso, Chiara and Grappein, Denise and Preziosi, Luigi and Scial{\`o}, Stefano},
  journal={Computers \& Mathematics with Applications},
  volume={151},
  pages={252--270},
  year={2023},
  publisher={Elsevier}
}

@article{berrone20263d,
  title={A 3D-1D virtual element method for modeling root water uptake},
  author={Berrone, Stefano and Ferraris, Stefano and Grappein, Denise and Teora, Gioana and Vicini, Fabio},
  journal={Computers \& Mathematics with Applications},
  volume={213},
  pages={22--52},
  year={2026},
  publisher={Elsevier}
}

@article{berrone2025effective,
  title={Effective polygonal mesh generation and refinement for VEM},
  author={Berrone, Stefano and Vicini, Fabio},
  journal={Mathematics and Computers in Simulation},
  volume={231},
  pages={239--258},
  year={2025},
  publisher={Elsevier}
}

@article{SORGENTE2022151,
title = {Polyhedral mesh quality indicator for the Virtual Element Method},
journal = {Computers \& Mathematics with Applications},
volume = {114},
pages = {151-160},
year = {2022},
issn = {0898-1221},
doi = {https://doi.org/10.1016/j.camwa.2022.03.042},
url = {https://www.sciencedirect.com/science/article/pii/S0898122122001390},
author = {T. Sorgente and S. Biasotti and G. Manzini and M. Spagnuolo}
}

@book{comp_top,
  author = {Edelsbrunner, Herbert and Harer, John},
  ee = {http://www.ams.org/bookstore-getitem/item=MBK-69},
  isbn = {978-0-8218-4925-5},
  pages = {I-XII, 1-241},
  publisher = {American Mathematical Society},
  title = {Computational Topology - an Introduction.},
  year = 2010
}

@article{voro,
author = {Du, Qiang and Faber, Vance and Gunzburger, Max},
title = {Centroidal Voronoi Tessellations: Applications and Algorithms},
journal = {SIAM Review},
volume = {41},
number = {4},
pages = {637-676},
year = {1999},
doi = {10.1137/S0036144599352836},
URL = {https://doi.org/10.1137/S0036144599352836},
eprint = {https://doi.org/10.1137/S0036144599352836}
}

@article{tetgen,
author = {Si, Hang},
title = {TetGen, a Delaunay-Based Quality Tetrahedral Mesh Generator},
year = {2015},
issue_date = {January 2015},
publisher = {Association for Computing Machinery},
address = {New York, NY, USA},
volume = {41},
number = {2},
issn = {0098-3500},
url = {https://doi.org/10.1145/2629697},
doi = {10.1145/2629697},
journal = {ACM Trans. Math. Softw.},
month = feb,
articleno = {11},
numpages = {36}
}

@Inbook{deBerg2008,
author="de Berg, Mark
and Cheong, Otfried
and van Kreveld, Marc
and Overmars, Mark",
title="Binary Space Partitions",
bookTitle="Computational Geometry: Algorithms and Applications",
year="2008",
publisher="Springer Berlin Heidelberg",
address="Berlin, Heidelberg",
pages="259--281",
isbn="978-3-540-77974-2",
doi="10.1007/978-3-540-77974-2_12",
url="https://doi.org/10.1007/978-3-540-77974-2_12"
}

@article{ANTONIETTI2022111531,
title = {Machine learning based refinement strategies for polyhedral grids with applications to virtual element and polyhedral discontinuous Galerkin methods},
journal = {Journal of Computational Physics},
volume = {469},
pages = {111531},
year = {2022},
issn = {0021-9991},
doi = {https://doi.org/10.1016/j.jcp.2022.111531},
url = {https://www.sciencedirect.com/science/article/pii/S0021999122005939},
author = {P.F. Antonietti and F. Dassi and E. Manuzzi}
}

@article{BERRONE2021103502,
title = {Refinement strategies for polygonal meshes applied to adaptive VEM discretization},
journal = {Finite Elements in Analysis and Design},
volume = {186},
pages = {103502},
year = {2021},
issn = {0168-874X},
doi = {https://doi.org/10.1016/j.finel.2020.103502},
url = {https://www.sciencedirect.com/science/article/pii/S0168874X20301827},
author = {Stefano Berrone and Andrea Borio and Alessandro D'Auria}
}

@article{KNUPP2003217,
title = {Algebraic mesh quality metrics for unstructured initial meshes},
journal = {Finite Elements in Analysis and Design},
volume = {39},
number = {3},
pages = {217-241},
year = {2003},
issn = {0168-874X},
doi = {https://doi.org/10.1016/S0168-874X(02)00070-7},
url = {https://www.sciencedirect.com/science/article/pii/S0168874X02000707},
author = {Patrick M. Knupp}
}

\end{document}